\documentclass[12pt]{article}
\usepackage{amsmath,amssymb,amscd,amsthm,latexsym,amsfonts}
\newtheorem{theorem}{Theorem}
\newtheorem{lemma}[theorem]{Lemma}
\newtheorem{proposition}[theorem]{Proposition}

\newtheorem{definition}[theorem]{Definition}

\newtheorem{remark}[theorem]{Remark}

\def\ep{\varepsilon}
\def\supp{\,{\rm Supp}\,}

\newcommand{\NN}{\mathbb{N}}
\newcommand{\RR}{\mathbb{R}}

\newcommand{\CC}{\mathbb{C}}

\def \reel{\RR}
\def \nat{ \NN}

\def\complex{\CC}

\renewcommand{\max}{\operatorname{max}}

\begin{document}
\author{\begin{tabular}{ccc}
Pascal Auscher \footnote{Universit\'e de Picardie-Jules Verne,   LAMFA,
CNRS FRE 2270, 33 rue Saint-Leu, F-80039 AMIENS Cedex 1,  E-mail:
auscher@u-picardie.fr}   & Emmanuel Russ\footnote{Facult\'e des
Sciences et Techniques de Saint-J\'er\^ome,
 Avenue Escadrille Normandie-Ni\'emen, F-13397 MARSEILLE Cedex 20
and LATP, CNRS UMR 6632. Email: emmanuel.russ@math.u-3mrs.fr}
\end{tabular}}
\title{Hardy spaces and divergence operators on strongly Lipschitz domains
of $\reel^n$}

\date{December 12, 2001}

\maketitle
\begin{abstract}
Let $\Omega$ be a strongly Lipschitz domain of $\reel^n$. Consider an
elliptic second order divergence operator
$L$ (including a boundary condition on $\partial\Omega$) and define a 
Hardy space by imposing the non-tangential maximal function of 
the  extension of a function $f$ via the Poisson semigroup for $L$ to
be in
$L^1$. Under suitable assumptions on $L$, we identify  this maximal
Hardy space with atomic Hardy spaces, namely with $H^1(\reel^n)$ if
$\Omega=\reel^n$,  
$H^{1}_{r}(\Omega)$ under the Dirichlet boundary condition, and
$H^{1}_{z}(\Omega)$ under the Neumann boundary condition. In
particular, we obtain a new proof of the atomic decomposition for 
$H^{1}_{z}(\Omega)$. A  version for local Hardy spaces is also given.
We also present an overview of the theory of Hardy spaces and BMO spaces
on Lipschitz domains with proofs.


MSC 2000 Classification numbers: 42B30, 42B25 

\end{abstract}
\vfill\break

\tableofcontents 

\vfill\break

\section{Introduction}

Hardy spaces on $\reel^n$, and especially $H^{1}(\reel^n)$, were
studied in great detail in the 60's and 70's.  A nice review on this is in
\cite{semmes}. 

  Originally
defined by means of Riesz transforms (see the seminal paper of Stein
and Weiss \cite{steinweiss}),  the usefulness of this space in
analysis as a substitute for
$L^1(\reel^n)$ comes from its many characterizations, beginning from the work 
of Fefferman-Stein  (see
\cite{feffstein}). Let
$\phi\in {\cal S}(\reel^n)$ be a function such that
$\int_{\reel^n} \phi(x)dx=1$. For all $t>0$, define $\phi_t(x)=t^{-n}
\phi(x/t)$. A locally integrable function $f$ on $\reel^n$ is said to be
in $H^1(\reel^n)$ if the vertical maximal function
\[
{\cal M}f(x)=\sup\limits_{t>0} \left\vert \phi_t \ast f(x)\right\vert
\]
belongs to $L^1(\reel^n)$. If it is the case, define
\[
\left\Vert f\right\Vert_{H^1(\reel^n)} = \left\Vert {\cal M}f\right\Vert_1.
\]
Recall that a function $f\in H^1(\reel^n)$ satisfies $\int_{\reel^n}
f(x)dx=0$.

Another equivalent definition of $H^{1}(\reel^n)$
 involves the non-tangential maximal function
associated with the Poisson semigroup (or the heat semigroup) generated by
$\Delta$, the Laplace operator on $\reel^n$. If $f\in
L^{1}_{loc}(\reel^n)$, the following are equivalent:
\begin{equation} \label{semigroup}
\begin{array}{l}
f\in H^{1}(\reel^n),\\
\\
\sup\limits_{\left\vert y-x\right\vert \leq t}
\left\vert e^{-t(-\Delta)^{1/2}}f(y)\right\vert \in L^{1}(\reel^n).
\end{array}
\end{equation}
See \cite{feffstein}, Theorem 11, p. 183.

The
atomic decomposition obtained by Coifman and Latter was a key step in
the theory (see
\cite{coif} when
$n=1$,
\cite{latter} when $n\geq 2$). A function $a$ on $\reel^n$ is
an $H^{1}(\reel^n)$-atom if it is supported in a cube $Q$, has mean-value
zero and
satisfies $\left\Vert a\right\Vert_{2}\leq {\left\vert
Q\right\vert^{-1/2}}$. Then,
$f \in H^{1}(\reel^n)$ if and only if $f=
\sum_{Q}
\lambda_{Q}a_{Q}$ where
the $a_Q$'s are $H^{1}(\reel^n)$-atoms  and the sequence  of complex
numbers $(\lambda_{Q})_{Q}$ is  in $l^{1}$. The norm
 $\left\Vert f\right\Vert_{H^{1}(\reel^n)}$ is
comparable with the infimum of
$\sum_{Q} \left\vert \lambda_{Q}\right\vert$ taken over all
such decompositions.

In  recent years, a quite complete theory of Hardy spaces on domains has
been developed  (\cite{jsw, miyachi, stkrch, dafni, triebel}.  The
Hardy  spaces are defined in terms of restrictions or support
conditions from
$H^1(\reel^n)$ or in terms of some ``grand' maximal function. For these
spaces, atomic decomposition  have
been obtained in particular on special  Lipschitz domains and bounded
Lipschitz domains of $\reel^n$. However, {\it is there a maximal
characterization using the Poisson semigroup?} More precisely, replace
in  (\ref{semigroup}),
$\reel^n$ by $\Omega$ and take for $\Delta$ the Laplacian with Dirichlet
or Neumann boundary condition. This defines two maximal Hardy spaces on
$\Omega$. One of the aims of the present paper is to identify each one
with one of the ``geometrical'' Hardy spaces mentioned above. It turns
out that the choice of boundary condition is meaningful  in the answer.
Roughly, the maximal space corresponding to the Dirichlet Laplacian is
$H^1_r(\Omega)$ and to the Neumann Laplacian 
$H^1_z(\Omega)$. In the  Dirichlet case, we shall use the existing
atomic decomposition of $H^1_r(\Omega)$. On the other hand, 
in the  Neumann case, we obtain in passing the atomic  decomposition
of $H^1_z(\Omega)$. We  also make
the statements valid for general strongly Lipschitz domains (which
include also exterior domains).

Another question we ask here is: {\it does the Laplacian play a specific
role?} In
other words, can it be replaced by an other second order elliptic operator? In
\cite{grenoble}, it was shown that $H^1(\RR)$ has a 
maximal characterization using the Poisson semigroup of 
elliptic operators. We give here an affirmative answer in higher
dimensions and on domains,  provided  the elliptic operator satisfies 
a technical condition.  For example any real elliptic
operator will do. This also emphasizes the prominent role of the
boundary condition in these questions.

The similar questions for local Hardy spaces have 
comparable answers.

Using the recent work of Dafni et al, \cite{dafni}, one can certainly
extend our results to $H^p$ and $h^p$ spaces for a range of $p$'s
smaller than 1.  We have not done so to keep the length of the paper
reasonable.

The plan of this paper is the following. First, we treat the case of
global Hardy spaces on strongly Lipschitz domains: we review their
definitions   and
recall their atomic decompositions (and clarify some points in the
literature).  We then introduce our maximal Hardy spaces and state the
main theorem. Next, we recall a few facts about
$BMO$ and duality.  Then we turn to proving some auxiliary results
involving square functions, Carleson measures and tent spaces before
proving the main theorem. We also give  proofs (sometimes new) of
classical atomic decomposition and of duality.   In a second part, we study
the corresponding theory for local Hardy spaces.  We also present
different maximal functions characterizing our maximal Hardy space. We
conclude with two appendices, one about kernel estimates and the other
about the elementary geometry of Lipschitz domains.

\paragraph{Acknowledgements} We are grateful to Alan McIntosh and  Zengjian
Lou for letting us include an argument of theirs and for discussions on this
project. The second author  thanks Philippe Tchamitchian for advice and
encouragement.

\section{Global Hardy  spaces  on strongly Lipschitz
domains}

In what follows,  it is understood without mention that $\Omega$ belongs
to the class of  strongly Lipschitz domains of
$\reel^n$, that is $\Omega$ is a proper open connected set in $\reel^n$
whose boundary is a finite union of parts of rotated graphs of Lipschitz
maps, at most one of these parts possibly  infinite. 

This class includes special Lipschitz domains, bounded Lipschitz
domains and exterior domains. Some facts about such domains are
presented in Appendix B.

Some statements may be valid for a restricted class and we shall
indicate this when it is the case.

\subsection{Hardy spaces: definitions}\label{sec:Hardydef}

Let us begin with defining  various Hardy spaces on a domain. Some
definitions will differ from (\cite{stkrch, dafni}). 
For the atomic spaces, we have privileged  
$L^2$ normalized atoms. We will not address the equivalent
definitions obtained by taking $L^p$ normalized atoms with $p>1$. 

\medskip

The first category is made up of  restrictions to $\Omega$ of certain
functions in $H^1(\reel^n)$.

\paragraph{Definition of $H^1_r(\Omega)$:} A function $f$ on $\Omega$
is said to be in $H^1_r(\Omega)$ if it is the restriction to $\Omega$ of a
function
$F\in H^1(\reel^n)$. If $f\in
H^1_r(\Omega)$, define
$\left\Vert f\right\Vert_{H^1_r(\Omega)}$ by
\[
\left\Vert f\right\Vert_{H^1_r(\Omega)}=\inf\left\Vert
F\right\Vert_{H^1(\reel^n)},
\]
the infimum being taken over all the functions $F\in H^1(\reel^n)$ such
that $\left.F\right\vert_{\Omega}=f$.

\paragraph{Definition of $H^1_z(\Omega)$:}A function $f$ on $\Omega$
belongs to $H^1_z(\Omega)$ if the function $F$ defined by
\[
F(x)=\left\{ 
\begin{array}{cc}
f(x) & \mbox{if }x\in \overline{\Omega},\\
\\
0 & \mbox{if }x\notin \Omega
\end{array}
\right.
\]
belongs to $H^1(\reel^n)$. When $f\in H^1_z(\Omega)$, its norm $\left\Vert
f\right\Vert_{H^1_z(\Omega)}$ is
$\left\Vert F\right\Vert_{H^1(\reel^n)}$. Note that it is a strict
subspace of $H^1_r(\Omega)$ (in particular, a function $f$ in
$H^1_z(\Omega)$ satisfies $\int_{\Omega}f(x)dx=0$,
whereas this may not happen for $f\in H^1_r(\Omega)$). This space is
nothing but the subspace of $H^1(\reel^n)$ of all functions supported in 
$\overline\Omega$.

\medskip

The second category of Hardy spaces on $\Omega$ consists of atomic
spaces. We list three such spaces.

\smallskip

\paragraph{Definition of type $(a)$ and $(b)$ cubes:} A cube $Q$ is said
to be a type $(a)$ cube [with respect to $\Omega$] if $4Q\subset
\Omega$,
a type $(b)$ cube if $2Q\subset \Omega$ and $4Q\cap
\partial \Omega\neq \emptyset$.

\paragraph{Definition of type $(a)$ and $(b)$ atoms:} A measurable
function $a$ on $\Omega$ is called a type $(a)$ atom if it is supported
in a type $(a)$ cube $Q$ with
\[
\int_Q a(x)dx=0\mbox{ and }\left\Vert a\right\Vert_2\leq \left\vert
Q\right\vert^{-1/2}.
\]
A measurable function $a$ on $\Omega$ is called a type $(b)$ atom if it is
supported in a type $(b)$ cube $Q$ with
\[
\left\Vert a\right\Vert_2\leq \left\vert Q\right\vert^{-1/2}.
\]
Note that a type $(b)$ atom is not supposed to have mean value zero. 

A more
speaking terminology would be {\it interior atoms} for type $(a)$ atoms and
{\it boundary atoms} for type
$(b)$ atoms. We have kept the terminology in the literature.

\paragraph{Definition of $H^{1}_{r,a}(\Omega)$:} A
function $f$ defined on $\Omega$ belongs to
$H^{1}_{r,a}(\Omega)$ if
\[
f=\sum\limits_{(a)} \lambda_{Q}a_{Q} + \sum\limits_{(b)} \mu_{Q}b_{Q}
\]
where the $a_Q$'s are type $(a)$ atoms, the $b_{Q}$'s are type
$(b)$ atoms and
$\sum_{(a)}\left\vert \lambda_Q\right\vert +
\sum_{(b)}\left\vert \mu_{Q}\right\vert <+\infty$.
Define $\left\Vert f\right\Vert_{H^1_{r,a}}$ as the infimum of
$\sum_{(a)}\left\vert \lambda_Q\right\vert +
\sum_{(b)}\left\vert \mu_{Q}\right\vert$ over
all such decompositions. 

\paragraph{Definition of $H^{1}_{z,a}(\Omega)$:}
A function $f$ defined on $\Omega$ belongs to
$H^{1}_{z,a}(\Omega)$ if
\[
f=\sum\limits_{(a)} \lambda_{Q}a_{Q}
\]
where the $a_Q$'s are type $(a)$ atoms and
$\sum_{(a)}\left\vert \lambda_Q\right\vert <+\infty$.
Define $\left\Vert f\right\Vert_{H^1_{z,a}}$ as the infimum of
$\sum_{(a)}\left\vert \lambda_Q\right\vert $ over
all such decompositions.

Note that this definition gives a smaller space than the atomic space
  considered in \cite{stkrch}: there, the $a_Q$'s are 
taken as $H^1_z(\Omega)$-atoms, that is   $H^1(\reel^n)$-atoms
supported  in cubes
contained in $\Omega$.
 We shall show 
that our definition coincides with theirs  (this fact is
 implicit in \cite{stkrch} when $\Omega$ is special 
Lipschitz or bounded as pointed out in  \cite{dafni}, p. 1612). An
immediate advantage of our definition of  $H^{1}_{z,a}(\Omega)$  is the 
strict containment of
$H^{1}_{z,a}(\Omega)$ in  $H^{1}_{r,a}(\Omega)$ and the evident role of the
boundary of $\Omega$. If $\Omega$ were  arbitrary, that definition could be
vacuous for application as $H^{1}_{z,a}(\Omega)$  could be too small.

\paragraph{Definition of $H^{1}_{CW}(\Omega)$:} Finally, since $\Omega$
is strongly Lipschitz, it is a space of  homogeneous type and one may
also  consider on $\Omega$ the Hardy  space of Coifman and Weiss as 
defined in \cite{coifw}, which will be denoted in the sequel by 
$H^{1}_{CW}(\Omega)$.
An $H^1_{CW}(\Omega)$-atom is a function $a$ supported in
$Q\cap\overline\Omega$, where $Q$ is a cube centered in $\Omega$
(but not necessarily included in $\Omega$) and satisfying
\[
\int a(x)dx=0\mbox{ and }\left\Vert a\right\Vert_{2} \leq 
{\left\vert Q \cap\Omega\right\vert^{-1/2}}.
\]
  If $\Omega$ has finite measure,  the
constant  function $\frac 1{\left\vert \Omega\right\vert}$
 is not an atom with  
our definition in opposition with that of \cite{coifw}.

A function $f$ is in $H^1_{CW}(\Omega)$ if $f$ can be written as 
\[
f=\sum\limits_{Q} \lambda_Q a_Q.
\]
where the $a_Q$'s are  $H^1_{CW}(\Omega)$-atoms and 
$\sum_{Q} |\lambda_Q| <\infty$. The norm is defined as usual.
This space is also
a special case of the Hardy space $H^1(F)$ considered in
\cite{jsw} on closed sets $F$ of
$\reel^n$ with the Markov property (here, $F=\overline\Omega)$.

\begin{theorem} \label{H1equalities}
\begin{itemize}
\item[$(a1)$] $H^1_{r}(\Omega) \subset H^1_{r,a}(\Omega)$.
\item[$(a2)$] $H^1_{r,a}(\Omega)= H^1_{r}(\Omega)$ provided
$^c\Omega$ is unbounded.
\item[$(b1)$]
$H^1_{z,a}(\Omega) = H^1_{CW}(\Omega)$.
\item[$(b2)$]
$H^1_{z,a}(\Omega) =  H^1_{z}(\Omega) $.
\end{itemize} 
\end{theorem} 

Each inclusion is here a continuous embedding between Banach spaces.
 In this paper, one finds a self-contained proof.  But, as this is not the
main object of our paper, let us  comment on this result now and postpone
proofs till later.

Assertions $(a1)$, $(a2)$  are known results  when
$\Omega$ is a special Lipschitz domain or is bounded \cite{miyachi, 
stkrch} and are simple to prove. We note that the restriction on
$\Omega$ in $(a2)$ is necessary as a counterexample will show
(See Section 2.8).

Concerning $(b1)$, the embedding 
$H^1_{z,a}(\Omega) \subset H^1_{CW}(\Omega) $ is straightforward. 
The converse can be   obtained
(but this is not straightforward) by typical arguments in harmonic
analysis on abstract homogeneous spaces combined with the geometry of the
boundary: maximal functions, Calder\'on-Zygmund decomposition and Whitney
coverings. However, we shall present a quite interesting argument due to Lou
and McIntosh
\cite{lm}, which uses more the differential structure of $\reel^n$ (See
Section 2.8).

That $H^1_{z,a}(\Omega) \subset H^1_{z}(\Omega)$ in $(b2)$ is a triviality.
The  remaining embedding 
$H^1_{z}(\Omega) \subset H^1_{z,a}(\Omega)$   is the deepest of all. 
It is proved by a constructive method
in \cite{stkrch} on special Lipschitz domains and on
bounded Lipschitz domains for the local Hardy spaces (see
Section \ref{local}).  In \cite{chang},  it is derived from
an extension theorem by Jones for $BMO$ \cite{Jones} and duality. 
Another argument is to use the  result that
$H^1_{z}(\Omega)=H^1(\overline\Omega)$ in \cite{jsw}, combined
with $H^1_{CW}(\Omega)=H^1(\overline\Omega)$ and  $(b1)$.

A
byproduct of our maximal spaces 
 defined below (Section \ref{sec:maxspaces}) is another proof of the
embedding $H^1_{z}(\Omega) \subset H^1_{CW}(\Omega)$. 
 
Up until Section 2.8, we assume knowledge of  Theorem \ref{H1equalities} but 
$H^1_{z}(\Omega) \subset H^1_{z,a}(\Omega)$ which is proved in Section 2.7.

\subsection {Maximal Hardy spaces and statement of the main
result}\label{sec:maxspaces}

We introduce a
third category of Hardy spaces on $\Omega$  defined via maximal
functions associated with second order elliptic operators in divergence
form.
We  briefly describe these operators, the most typical being the
Laplacian with appropriate boundary condition. If
$\Omega=\reel^n$ or if 
$\Omega$ is a strongly Lipschitz domain of $\reel^n$, we will denote by $W^{1,2}(\Omega)$
the usual Sobolev space on $\Omega$
equipped with the norm $\left(\left\Vert f\right\Vert_{2}^2+ \left\Vert \nabla
f\right\Vert_{2}^2\right)^{1/2}$,
whereas $W^{1,2}_0(\Omega)$ stands for the closure of
$C^{\infty}_0(\Omega)$ in $W^{1,2}(\Omega)$.

If $A:\reel^n \rightarrow M_{n}({\complex})$ is a measurable function,
define
\[
\left\Vert A\right\Vert_{\infty}=\sup\limits_{x\in \reel^n,\
\left\vert \xi\right\vert=\left\vert \eta\right\vert=1}\left\vert \langle
A(x)\xi,\eta\rangle\right\vert.
\]
Here and subsequently in the paper, the notation sup is used for
esssup.
For all $\delta>0$, denote by ${\cal A}(\delta)$ the class of all
measurable functions $A:\reel^n \rightarrow M_{n}({\complex})$
satisfying, for all $x,\xi\in \reel^n$:
\[
\left\Vert A\right\Vert_{\infty} \leq \delta^{-1}\mbox{ and Re
}\langle A(x)\xi,\xi\rangle \geq \delta \left\vert \xi\right\vert^{2}.
\]
Denote by ${\cal A}$ the union of all ${\cal A}(\delta)$ for $\delta>0$.

When $A\in {\cal A}$ and $V$ is a closed subspace of $W^{1,2}(\Omega)$
containing $W^{1,2}_0(\Omega)$, denote by $L$ the maximal-accretive
operator on $L^2(\Omega)$ with largest domain ${\cal D}(L)\subset V$ such
that
\begin{equation} \label{accretive}
\langle Lf,g\rangle=\int_{\Omega} A\nabla f.\overline{\nabla g},\ \forall
f\in {\cal D}(L),\forall g\in V.
\end{equation}
We will write $L=(A,\Omega,V)$. Say that $L$ satisfies the Dirichlet boundary condition (DBC) when
$V=W^{1,2}_0(\Omega)$, the Neumann boundary
condition (NBC) when $V=W^{1,2}(\Omega)$.

We turn to the definition of maximal Hardy spaces associated with such operators.
Let $L=(A,\Omega,V)$ be  as
above. This operator has a unique maximal accretive square root
$L^{1/2}$ so that $-L^{1/2}$ is the generator of an
$L^2(\Omega)$-contracting semigroup $P_t=e^{-tL^{1/2}}, t>0$, the
Poisson semigroup for $L$. We will need that $P_t$  also acts on
$L^1(\Omega)$. Let us then introduce a technical condition on $L$.  

\begin{definition} For $0<\tau\le +\infty$, we call   
 $(G_{\tau})$ the conjunction of (\ref{Gauss}) and
(\ref{Holder}) below:
The kernel of $e^{-tL}$, denoted by $K_t(x,y)$, is a measurable 
function on $\Omega\times\Omega$ and there exist $C,\alpha>0$ such 
that, for all $0<t<\tau$ and almost every $x,y\in \Omega$,
\begin{equation} \label{Gauss}
\left\vert K_t(x,y)\right\vert \leq \frac {C}{t^{n/2}}e^{-\alpha
\frac{\left\vert x-y\right\vert^2}t}.
\end{equation} 

For all $x\in \Omega$ and all $0<t<\tau$, the 
function $y\mapsto K_t(x,y)$ is H\"older continuous in $\Omega$
and there exist $C,\mu\in ]0,1]$ such that,
for all $0<t<\tau$ and all $x,y,y^{\prime}\in \Omega,$
\begin{equation} \label{Holder}
\left\vert K_t(x,y)-K_t(x,y^{\prime})\right\vert \leq \frac {C}{t^{n/2}}
\frac{\left\vert
y-y^{\prime}\right\vert^{\mu}}{t^{\mu/2}}.
\end{equation}
\end{definition}

When $\tau$ is finite, we set $\tau=1$ without loss of generality.

For those readers only interested in the Laplacian or
 real symmetric operators (under BDC or NBC),
this condition  is always satisfied on $\reel^n$ or on Lipschitz domains
with $\tau=\infty$ except under NBC  with $\Omega$ bounded for which we
 have $\tau$ finite.

\begin{lemma}\label{Poissonbounds} When $(G_{\infty})$ holds, the
Poisson kernel of $L$, i.e. the kernel $p_t(x,y)$ of $P_t$ satisfies
\begin{equation} \label{Poisson}
\left\vert p_t(x,y)\right\vert \leq  \frac{{C}{t}}{( t+|x-y|)^{n+1}}
\end{equation}
and 
\begin{equation} \label{PoissonHolder}
\left\vert p_t(x,y)-p_t(x,y^{\prime})\right\vert \leq \frac {C}{t^{n}}
\frac{\left\vert
y-y^{\prime}\right\vert^{\mu}}{t^{\mu}}
\end{equation}
for all $t\in ]0,\infty[$, for some $C>0$ and $\mu\in ]0,1[$.
\end{lemma}

This follows from  the subordination
formula (see \cite{topics}):
 \begin{equation} \label{subor}
 p_t(x,y)=\frac 1{\sqrt{\pi}}\int_0^{+\infty}
 K_{\frac{t^2}{4u}}(x,y)e^{-u}u^{-1/2}du.
 \end{equation}
\rule{2mm}{2mm}

  If $(G_{\infty})$ holds and
$f\in L^1_{loc}(\Omega)$ so that $y\mapsto |y|^{-n-1}f(y) \in
L^1(\Omega)$, define, for all
$x\in
\Omega$,
\[
f^*_L(x)=\sup\limits_{y \in \Omega, t>0,\ \left\vert y-x\right\vert<t}
\left\vert P_tf(y)\right\vert.
\]
 Say that $f\in H^1_{max,L}(\Omega)$ if $f^*_L\in L^1(\Omega)$
and define
\[
\left\Vert f\right\Vert_{H^1_{max,L}(\Omega)}=\left\Vert
f^*_L\right\Vert_{L^1(\Omega)}.
\]
Note that $H^1_{max,L}(\Omega)$ depends, in
particular, on the boundary condition. Since $P_t$ tends to the identity
strongly in $L^1(\Omega)$ we see that $H^1_{max,L}(\Omega) \subset
L^1(\Omega)$.

One of the aims of  this paper is to identify this maximal space. 
Our  result is the following:
\begin{theorem} \label{comparison} 
Let $\Omega=\reel^n$ or $\Omega$ be a strongly Lipschitz
domain of $\reel^n$, and
$L=(A,\Omega,V)$ satisfying $(G_{\infty})$.
\begin{itemize}
\item[$(a)$]
If $\Omega=\reel^n$, one has $H^{1}(\reel^n)=H^{1}_{max,L}(\reel^n)$.
\item[$(b)$]
If $\, ^c \Omega$ is unbounded 
and 
$L$ satisfies the Dirichlet boundary condition, then one has
$H^1_{r,a}(\Omega)=H^1_{max,L}(\Omega)= H^1_{r}(\Omega)$.
\item[$(c)$]
If $\Omega$ is  unbounded  and $L$ satisfies the Neumann
boundary  condition, then one has
$H^1_{z,a}(\Omega)=H^1_{max,L}(\Omega)=H^1_z(\Omega)$.
\end{itemize}
\end{theorem}

If we assume $^c \Omega$ is
bounded in $(b)$ then $H^1_{r,a}(\Omega)\subset H^1_{max,L}(\Omega)$. We
have not succeeded in proving the converse. Assuming $\Omega$ unbounded in
$(c)$ is no restriction as the condition $(G_{\infty})$ in never
satisfied under NBC and 
$\Omega$ bounded.   Note that 
$(c)$ contains the equality between Hardy spaces alluded to in Section
\ref{sec:Hardydef} in the case where
$\Omega$ in unbounded.  The case where
$\Omega$ is bounded will be addressed via local spaces
in Section 2.1.   We turn to  some
intermediate results and begin with discussing about $BMO$ spaces.

\subsection{$BMO$ spaces}\label{BMO}  

\paragraph{Definition of $BMO(\reel^n)$:} A locally square-integrable
function
$f$ on
$\reel^n$ is said to be in
$BMO(\reel^n)$ if
\[
\left\Vert f\right\Vert_{BMO(\reel^n)}^2=\sup\limits_Q \frac 1{\left\vert
Q\right\vert} \int_Q \left\vert f(x)-f_Q\right\vert^2 dx<+\infty
\]
where the supremum is taken over all the cubes $Q\in \reel^n$ with sides
parallel to the axes. Here, $f_E=\frac 1{\left\vert E\right\vert} \int_E
f(x)dx$ is the mean of $f$ over $E$ and $|E|$ is the Lebesgue measure
of $E$.

\paragraph{Definition of $VMO(\reel^n)$:} Define $VMO(\reel^n)$
as  the closure
of $C_c(\reel^n)$ (the space of continuous functions on
$\reel^n$ with compact support) in
$BMO(\reel^n)$.  This $VMO$ space is the one in the sense of Coifman
and Weiss \cite{coifw} and is different from the one considered by
Sarason in \cite{sarason}, which is the closure of the space of
all uniformly continuous $BMO$-functions on
$\reel^n$. See the recent work of G. Bourdaud for clarifications
\cite{bourdaud}. It is well-known that $BMO(\reel^n)$ is the dual of
$H^1(\reel^n)$, the latter being the dual of $VMO(\reel^n)$ 
\cite{feffstein, coifw}.

\medskip

We next introduce the first category of $BMO$-spaces on $\Omega$.

\paragraph{Definition of $BMO_z(\Omega)$:} The space $BMO_z(\Omega)$
is defined as being the space of all functions
in $BMO(\reel^n)$ supported in $\overline{\Omega}$, equipped with the norm
$\left\Vert f\right\Vert_{BMO_z(\Omega)}=\left\Vert
f\right\Vert_{BMO(\reel^n)}$.

\paragraph{Definition of $VMO_z(\Omega)$:} We define $VMO_z(\Omega)$
is the closure of $C_c(\Omega)$, the space of continuous
functions  with support in $\Omega$, in $BMO_z(\Omega)$.

\paragraph{Definition of $BMO_r(\Omega)$:} The space $BMO_r(\Omega)$
is defined as being the space of all restrictions to $\Omega$ of functions 
in $BMO(\reel^n)$. If $f \in BMO_r(\Omega)$ define
$\left\Vert f\right\Vert_{BMO_r(\Omega)}$ by
\[
\left\Vert f\right\Vert_{BMO_r(\Omega)}=\inf\left\Vert
F\right\Vert_{BMO(\reel^n)},
\]
the infimum being taken over all the functions $F\in BMO(\reel^n)$ such
that $\left.F\right\vert_{\Omega}=f$.
\medskip

Next, we turn to the second category of $BMO$-spaces, defined in terms of mean
square oscillation.

\paragraph{Definition of $BMO_{z,a}(\Omega)$:}
A locally square-integrable function $f$ on $\Omega$ is in 
$BMO_{z,a}(\Omega)$ if
\[
\left\Vert \phi\right\Vert_{BMO_{z,a}(\Omega)}^2=\sup
\left(\sup\limits_{(a)}
\frac 1{\left\vert Q\right\vert} 
\int_{Q} \left\vert \phi(x)-\phi_{Q}\right\vert^2 dx 
, 
\sup\limits_{(b)}  \frac
1{\left\vert Q\right\vert} 
\int_{Q}\left\vert\phi\right\vert^2\right) <+\infty,
\]
where $\sup\limits_{(a)}$ (resp. $\sup\limits_{(b)}$) means 
that the supremum is taken over all 
  type $(a)$ cubes (resp. all 
 type $(b)$ cubes).

\paragraph{Definition of $BMO_{r,a}(\Omega)$:}
A locally square-integrable function $f$ on $\Omega$ is in 
$BMO_{r,a}(\Omega)$ if
\[
\left\Vert \phi\right\Vert_{BMO_{r,a}(\Omega)}^2=
\sup\limits_{(a)}
\frac 1{\left\vert Q\right\vert} 
\int_{Q} \left\vert \phi(x)-\phi_{Q}\right\vert^2 dx 
<+\infty.
\] 
\smallskip

\paragraph{Definition of $BMO_{CW}(\Omega)$:}
A locally square-integrable function $\phi$  on $\Omega$ is in
$BMO_{CW}(\Omega)$ if
\[
\left\Vert \phi\right\Vert_{BMO_{CW}(\Omega)}^2 = 
\sup \frac 1{\left\vert Q\cap\Omega\right\vert} \int_{Q\cap\Omega}
\left\vert \phi(x)-\phi_{Q\cap\Omega}\right\vert^2 dx<+\infty,
\]
where the supremum is taken over all cubes centered in
$\Omega$. This is the space defined in \cite{coifw}. A sligth
variation is that  the indicator function of $\Omega$
is not in $BMO_{CW}(\Omega)$ when  $\Omega$ has finite measure.

\paragraph{Definition of $VMO_{CW}(\Omega)$:} The space
$VMO_{CW}(\Omega)$ is the closure of
$C_c(\Omega)$ in $BMO_{CW}(\Omega)$.
\smallskip

Note that $BMO_{r,a}(\Omega)$, $BMO_{r}(\Omega)$ and 
${BMO_{CW}(\Omega)}$ are defined modulo constants. We ignore
this well-understood issue.

Let us mention the duality results. 

\begin{theorem} \label{thdual} 
\begin{itemize}
\item[$(a)$]
The dual of $H^1_{r,a}(\Omega)$ is $BMO_{z,a}(\Omega)$.
\item[$(b)$]
The dual of $H^1_{r}(\Omega)$ is $BMO_{z}(\Omega)$, the dual of 
$VMO_{z}(\Omega)$ is $H^{1}_{r}(\Omega)$.
\item[$(c)$]
The dual of $H^{1}_{CW}(\Omega)$ is $BMO_{CW}(\Omega)$, the dual of
$VMO_{CW}(\Omega)$ is $H^{1}_{CW}(\Omega)$. 
\item[$(d)$]
The dual of $H^1_z(\Omega)$ is $BMO_r(\Omega)$.
\item[$(e)$]
The dual of $H^1_{z,a}(\Omega)$ is $BMO_{r,a}(\Omega)$.
\end{itemize}
\end{theorem}

The result corresponding to Theorem \ref{H1equalities} for $BMO$-spaces
is the following

\begin{theorem} \label{BMOequalities}
\begin{itemize}
\item[$(a1)$]
$BMO_{z,a}(\Omega) \subset BMO_{z}(\Omega)$. 
\item[$(a2)$]
$BMO_{z,a}(\Omega)= BMO_{z}(\Omega)$ 
provided $^c\Omega$ is unbounded.  
\item[$(b1)$]
$BMO_{CW}(\Omega)=  BMO_{r,a}(\Omega) $.
\item[$(b2)$]
$BMO_{r}(\Omega) =   BMO_{r,a}(\Omega) $.
\end{itemize}\smallskip
\end{theorem}

Again, let us just comment on these results of which we shall only need 
$(b)$ and $(c)$ of Theorem \ref{thdual}, which will be proved in Section
\ref{duality}.

Concerning Theorem \ref{thdual}, $(c)$ is already known \cite{coifw}
(our change in the definition does not induce any modification in
thte proof) and all the other statements but
$(e)$ are not so deep. It is easy to show that $(H^1_{z,a}(\Omega))'$
contains
$BMO_{r,a}(\Omega)$ but the converse is harder (in particular, we
think that the argument proposed in \cite{chang}, Theorem 2.1, for
bounded domains and local  spaces has a gap). 

Assume for the moment all the
above is proved and let us argue for Theorem \ref{BMOequalities}.
First, $(a1)$ and $(a2)$  follow by  duality (easy direct proofs are
also possible). Next and we have
$BMO_{r}(\Omega) \subset BMO_{CW}(\Omega)\subset  
BMO_{r,a}(\Omega)$ (duality or direct proof). Using Theorem
\ref{H1equalities}, $(b1)$, and the already observed inclusion in
Theorem
\ref{thdual}, $(e)$, we have
$$
BMO_{CW}(\Omega) \subset BMO_{r,a}(\Omega) \subset (H^1_{z,a}(\Omega))'
= (H^1_{CW}(\Omega))' = BMO_{CW}(\Omega).$$
Hence,  $ BMO_{CW}(\Omega)=BMO_{r,a}(\Omega)$ and $
BMO_{r,a}(\Omega)$ is the dual of
$H^1_{z,a}(\Omega)$. 

This completes the proof of Theorem \ref{thdual}, $(e)$, and of 
Theorem \ref{BMOequalities}, $(b1)$.

 The remaining embedding in Theorem \ref{BMOequalities}, $(b2)$,
namely, 
$BMO_{r,a}(\Omega) \subset BMO_{r}(\Omega)$, given  duality and the
existing embeddings, is equivalent to  the embedding $H^1_z(\Omega)
\subset H^1_{z,a}(\Omega)$.

A result by Jones \cite{Jones}
characterizes the domains having an extension property for $BMO$.
Lipschitz domains fall in that class.   The embedding
$BMO_{r,a}(\Omega) \subset BMO_{r}(\Omega)$ is a slightly stronger 
extension property     since the definition of 
$BMO_{r,a}(\Omega)$ requires bounded mean oscillation only on cubes of
type
$(a)$ while Jones assumes  bounded mean
oscillation on all cubes inside  $\Omega$.

\subsection{Area integrals and maximal functions} \label{area}

Let $\Omega=\reel^n$ or $\Omega$ be a strongly Lipschitz domain of 
$\reel^n$.
Consider $L=(A,\Omega,V)$ with Dirichlet or Neumann boundary
condition.    Define for
$x\in
\Omega$, 
\[
S_{\alpha}f(x)=\left(\int_{\Gamma_{\alpha}(x)}
t^{1-n}\left\vert
\overline{\nabla}P_tf(y)\right\vert^2 dydt\right)^{1/2},
\]
and
\[
S_{\alpha}^{\ep,R}f(x)=\left(\int_{\Gamma_{\alpha}^{\ep,R}(x)}
t^{1-n}\left\vert
\overline{\nabla}P_tf(y)\right\vert^2 dydt\right)^{1/2},
\]
with $\overline{\nabla}u= (\nabla u, \partial_{t}u) $,
$\left\vert
\overline{\nabla}u\right\vert^2=
 \left\vert \nabla
u \right\vert^2+ \left\vert \partial_{t}u \right\vert^2$, 
$P_{t}=e^{-tL^{1/2}}$ and where
$\Gamma_{\alpha}(x)$ and $\Gamma_{\alpha}^{\ep,R}(x)$  are
the respectively the cones and the truncated cones defined by
\[
\Gamma_{\alpha}(x)=\left\{(y,t)\in \Omega\times
\left]0,+\infty\right[;\left\vert y-x\right\vert 
<\alpha t\right\}
\]
and 
\[
\Gamma_{\alpha}^{\ep,R}(x)=\left\{(y,t)\in \Omega\times
\left]\ep,R\right[;\left\vert y-x\right\vert <\alpha
t\right\},
\]
for   $\alpha>0$, $0<\ep<R<+\infty$. We can also write
$S_\alpha=S_\alpha^{0,\infty}$.  
Here, $|\ |$ is the sup norm on $\reel^n$ (for which the
balls are cubes with sides parallel to the axes).

\begin{lemma}\label{technicallemma} Assume  $\alpha <1$. Then,
one has for
$f\in L^2(\Omega)$,
$$
S_{\alpha}^{\ep,R}f(x) \le C (1+|\ln(R/\ep)|) f^*_L(x)
$$
for some constant depending on $\alpha$.
\end{lemma}

The truncated square function is
well-defined for 
$f \in L^2(\Omega)$ since $\overline\nabla P_t$ is bounded on
$L^2(\Omega	)$. Let us also recall that $u_t(y)=P_tf(y)$
satisfies 
the elliptic equation $\overline\nabla \cdot B\overline\nabla
u_t(y)=0$ (in the weak sense on $ \Omega \times ]0,\infty[$)
where $B$ is the $(n+1)\times (n+1)$ block diagonal matrix
with components $A$ and 1. Moreover, we have prescribed Dirichlet or
Neumann data on the lateral boundary $\partial\Omega \times
]0,\infty[.$ Hence, we have interior and boundary Caccioppoli
inequalities (see \cite{autchdom}): for some $\rho>0$ and $C$
depending on $\Omega$ and ellipticity, 
$$\int_{E} |\overline{\nabla}u_t(y)|^2 \, dydt \le C
r^{-2} \int_{\tilde E} |u_t(y)|^2 \, dydt
 $$
for all sets $E=B((z,\tau),r) \cap  (\Omega \times ]0,\infty[)$
with $\tilde E= B((z,\tau), 2r) \cap  (\Omega \times ]0,\infty[)$
provided $x \in \Omega$, $\tau>0$ and 
$r \le \inf (\rho, \tau)/4$. Here $B((z,\tau),r)$ is the open ball
defined by $\sup(|z-y|,|\tau - t|) < r$.

For $(z,\tau) \in \Gamma_{\alpha}^{\ep,R}(x)$, let $E_{(z,\tau)}=
B((z,\tau),r) \cap  (\Omega \times ]0,\infty[)$ with $r=\delta \inf(
\tau, \rho)$  . Here $\delta$ is some small number. By Besicovitch
covering argument,    pick a  subcollection 
$E_j=E_{(z_j,\tau_j)}$ covering $\Gamma_{\alpha}^{\ep,R}(x)$ and having
bounded overlap. Remark that $(y,t) \in E_j$ implies $t \sim d_j$,
the distance from $E_j$ to the bottom boundary $\Omega\times \{0\}$.
Remark also that if $\delta$ is small enough, $(y,t) \in \tilde E_j$
implies 
$(y,t) \in  \Gamma_{1}(x)$, hence $|u_t(y)| \le f^*_L(x)$.
Thus we obtain from the bounded overlap and Caccioppoli's inequality,   
$$
S_{\alpha}^{\ep,R}f(x)^2 \le C \sum_j d_j^{1-n} r_j^{-2} |\tilde E_j|
f^*_L(x)^2.$$
Observe that  $|\tilde E_j| \le C |E_j|$ and so that the bounded
overlap of the $E_j$'s again easily yields by inspection,
$$\sum_j d_j^{1-n} r_j^{-2} |\tilde E_j| \le 
C(1+|\ln(R/\ep)|).$$
\rule{2mm}{2mm}

\begin{proposition} \label{areaintegral} Assume that 
$(G_{\infty})$ holds.
There exists $C>0$
such that, for all $f\in H^1_{max,L}(\Omega)$,
$\left\Vert S_1f\right\Vert_{1}
\leq C\left\Vert f\right\Vert_{H^1_{max,L}}$.

\end{proposition}

The proof 
follows  ideas from \cite{feffstein}, Theorem 8, p. 161 and 
\cite{tent}, Section 6, see also \cite{grenoble}, Lemme II.10.
 It relies on 
a ``good $\lambda$'' inequality. We need though  variants of the
truncated square functions in order to compensate the lack of
pointwise regularity. Set
\[
\widetilde
S_{\alpha}^{\ep,R}f(x)=\left(\int_1^2
\int_{\Gamma_{\alpha/a}^{a\ep,aR}(x)}
t^{1-n}\left\vert
\overline{\nabla}P_tf(y)\right\vert^2 dydtda\right)^{1/2}, \quad  x
\in \Omega.
\]
Fairly elementary arguments show that
$$
S_{\alpha}^{2\ep,R}f \le \widetilde
S_{\alpha}^{\ep,R}f \le S_{2\alpha}^{\ep,2R}f.
$$
We shall prove

\begin{lemma} \label{distri}
There exists  $c>0$ such that, for all
$0<\gamma\leq 1$, all $\lambda>0$, all $0<\ep<R<\infty$ and
all
$f
\in H^1_{max,L} \cap L^2(\Omega)$,
\[
\left\vert \left\{x\in \Omega ;\widetilde
S_{1/20}^{\ep,R}f(x)>2\lambda,\ f^{\ast}(x)\leq \gamma
\lambda\right\} \right\vert \leq c\gamma^{2} \left\vert
\left\{x\in \Omega ;\widetilde
S_{1/2}^{\ep,R}f(x)>\lambda\right\} \right\vert.
\]
\end{lemma}

We will also  use  the comparability of
the square functions. See \cite{tent}, Proposition 4, p. 309.

\begin{lemma} \label{aperture}
For $\alpha,\beta>0$, $0\le \ep<R\le +\infty$, one has
\[
\left\Vert S_{\alpha}^{\ep,R}f\right\Vert_{1} \sim \left\Vert
S_{\beta}^{\ep,R}f\right\Vert_{1},
\]
where the implicit constants do not depend on $f, \ep, R$.
\end{lemma}

Let us deduce Proposition \ref{areaintegral}. Assume first that
$f\in H^1_{max,L}
\cap L^2(\Omega)$.  As a consequence of Lemma \ref{distri}, by
integrating both sides with respect to $\lambda$, one obtains
\[
\left\Vert \widetilde
S_{1/20}^{\ep,R}f\right\Vert_{1} \leq \gamma^{-1}\left\Vert
f^{\ast}_L\right\Vert_{1}+c\gamma^{2}\left\Vert \widetilde
S_{1/2}^{\ep,R}f\right\Vert_{1}.
\]
Thanks to Lemma \ref{aperture} and the comparisons between
the square functions and their variants, one has
$$\left\Vert S_{1}^{\ep,R}f\right\Vert_{1} \le C\left\Vert
\widetilde S_{1/20}^{\ep/2,R}f\right\Vert_{1}
$$
and by Lemma  \ref{technicallemma}
$$\left\Vert \widetilde
S_{1/2}^{\ep/2,R}f\right\Vert_{1} \le \left\Vert 
S_{1}^{\ep/2,2R}f\right\Vert_{1} \le \left\Vert
 S_{1}^{\ep/2,\ep}f\right\Vert_{1}
+ \left\Vert 
S_{1}^{\ep,R}f\right\Vert_{1}
+ \left\Vert 
S_{1}^{R,2R}f\right\Vert_{1} \le 
\left\Vert
 S_{1}^{\ep,R}f\right\Vert_{1} + C\left\Vert
f^{\ast}_L\right\Vert_{1}.
$$
Hence, by choosing $\gamma$ appropriately and using the a
priori knowledge that $\left\Vert
 S_{1}^{\ep,R}f\right\Vert_{1} <+\infty$
 one obtains 
$$
\left\Vert
 S_{1}^{\ep,R}f\right\Vert_{1} \le C\left\Vert
f^{\ast}_L\right\Vert_{1}.
$$
By letting $\ep \downarrow 0$ and $R\uparrow +\infty$, the
conclusion of  Proposition
\ref{areaintegral} in the case  $f\in L^{2}$
follows.

To  complete the proof of Proposition \ref{areaintegral},
we have to relax the assumption $f\in L^2(\Omega)$. But, if
$f^*_L\in L^1$, then $f\in L^1(\Omega)$ and together with the
kernel estimates on the kernel of $P_t$,  one has
$f_{\varepsilon}\in L^2(\Omega)$ for all
$\varepsilon>0$, where
$f_{\varepsilon}(x)=P_{\varepsilon}f(x)$.
It follows that $\left\Vert Sf_{\varepsilon}\right\Vert_1 \leq C\left\Vert
(f_{\varepsilon})^*_L\right\Vert_{1}\leq \left\Vert
f^{\ast}_L\right\Vert_{1}$. Letting $\varepsilon \downarrow
0$, one obtains
$\left\Vert Sf\right\Vert_1 \leq C\left\Vert
f^*_L\right\Vert_{1}$ by monotone convergence.
\rule{2mm}{2mm}

We turn to the proof of Lemma \ref{distri}. In the next
argument, $\ep,R,\lambda$ are fixed. Also $f \in H^1_{max,L}
\cap L^2(\Omega)$. Define
$O=\left\{x\in
\Omega  ;\widetilde
S_{1/2}^{\ep,R}f(x)>\lambda\right\}$. We may assume that
$O\neq \Omega$. Let $O=\bigcup\limits_{k}Q_{k}$ be a Whitney
decomposition of $O$ (with respect to $\Omega$) by dyadic
cubes (of $\reel^n$), so that, for all $k$,  
$2Q_{k}\subset O\subset
\Omega$, but $4Q_k$ intersects $\Omega\setminus O$. Since
 $\left\{\widetilde
S_{1/20}^{\ep,R}f>2\lambda\right\}\subset
\left\{\widetilde
S_{1/2}^{\ep,R}f>\lambda\right\}$, it is enough to show
that
\[
\left\vert \left\{x\in Q_{k};\ \widetilde
S_{1/20}^{\ep,R}f(x)>2\lambda,\
f^{\ast}_L(x)\leq \gamma
\lambda\right\} \right\vert \leq c\gamma^{2} \left\vert
Q_{k}\right\vert.
\]
From now on, fix $k$ and denote by $l$ the side length of 
$Q_{k}$.

If $x\in Q_{k}$,
\[
\widetilde
S_{1/20}^{\sup(10l,\ep),R}f(x)\leq  \lambda.
\]
Indeed, pick $x_k\in 4Q_{k}$ with $x_k \notin O$. If 
$\left\vert y-x\right\vert < 
\frac t{20} $ and 
with $t\geq  \sup(10l,\ep)$, then one has 
$\left\vert
x_{k}-y\right\vert < \frac t{20} + 4l\le \frac t2 
 $. Hence $ \widetilde
S_{1/20}^{\sup(10l,\ep),R}f(x) \le \widetilde
S_{1/2}^{\sup(10l,\ep),R}f(x_k) \le \lambda.$

If $\ep \ge 10l$, we are done. Otherwise, 
using, $\widetilde
S_{1/20}^{\ep,R}f(x) \le \widetilde
S_{1/20}^{\ep,10l}f(x)+\widetilde
S_{1/20}^{10l,R}f(x)$, it remains to show that
$$
\left\vert \left\{x\in Q_{k} \cap F;\ g(x) >\lambda\right\}
\right\vert
\leq c\gamma^{2}
\left\vert Q_{k}\right\vert
$$
where
$$
g(x)= \widetilde
S_{1/20}^{\ep,10l}f(x)
$$
and 
\[
F=\left\{x\in \Omega;f^{\ast}_L(x)\leq \gamma\lambda\right\}.
\]
By Tchebytchev's inequality, this follows from 
$$
\int_{Q_k\cap F} g^2 \le c\gamma^{2} \lambda^2
\left\vert Q_{k}\right\vert.
$$
We note that the condition
$(G_{\infty})$ implies that $F$ is a closed set of
$\Omega$. 

If $5l\le \ep$, then the  argument using Caccioppoli's
inequality shows that
$$
\int_{Q_k\cap F} g^2 \le c\int_{Q_k\cap F} (f^*_L)^2 
\le c\gamma^2\lambda^2 |Q_k\cap F|.
$$

Assume  from now on that $\ep<5l$.
By geometric considerations,
\[
\int_{Q_{k}\cap F} g(x)^{2}dx \leq  
c\int_{1}^{2}\int_{{\cal E}_{a}} t\left\vert
\overline{\nabla}u_t(y)\right\vert^2 dydtda,
\]
where 
\[
{\cal E}_{a}=\left\{(y,t) \in \Omega \times ]\ep a, 10l a[; 
a \psi(y) <t\right\}
\]
with $\psi(y)$ the Lipschitz function equal to $20
\,{\rm dist}\,(y,Q_k\cap F)$. Recall also that
$u_t(y)=P_tf(y)$.

Observe that ${\cal E}_{a} = \{(y,at); (y,t) \in {\cal
E}_{1}\}.$ Define $E=\{y; (y,t) \in {\cal E}_{1}\}$: this is
an open set   in $\Omega$. For  a
connected component
$C$ of $E$, we let ${\cal C}_{a}= \{(y,t) \in {\cal E}_{a}; y
\in C\}$. It suffices to show that
$$
\int_{1}^{2}\int_{{\cal C}_{a}} t\left\vert
\overline{\nabla}u_t(y)\right\vert^2 dydtda \le
c\gamma^2\lambda^2 |C|.
$$
Indeed, summing over all connected components of $E$, we get
$$
\int_{1}^{2}\int_{{\cal E}_{a}} t\left\vert
\overline{\nabla}u_t(y)\right\vert^2 dydtda \le
c\gamma^2\lambda^2 |E|,
$$
and it remains to observe that  $E \subset 2Q_k$. Indeed, if
$y \in E$, there is a point $(y,t) $ above  contained ${\cal
E}_1$, hence there exists $x \in Q_k\cap F$ such that
$|y-x| < \frac t{20}$. Since $t<10l$, we have $|y-x| <\frac l2$ and
the desired inclusion follows.

We next fix a connected component $C$ of $E$. 
Consider
$a\in \left] 1, 2\right[$ and note that ${\cal C}_{a}$ is 
connected and has Lipschitz boundary. The ellipticity condition for
$A$ shows that
\[
\int_{{\cal C}_{a}} t\left\vert
\overline{\nabla}u_t(y)\right\vert^2 dydt  \leq  C\mbox{ Re
}\int_{{\cal C}_{a}}   tB\overline\nabla u_t(y)
\overline{\overline\nabla u_t(y)}\,dydt= C\mbox{ Re 
}I_{a},
\]
where $B$ is the $(n+1)\times (n+1)$ block diagonal matrix
with components $A$ and 1. The function $u_t(y)$ satisfies 
the equation $\overline\nabla \cdot B\overline\nabla
u_t(y)=0$ (in the weak sense on $ \Omega \times ]0,\infty[$) so that 
we wish to integrate by parts. 

To do so let us make some
observations. We claim that for $(y,t) \in \overline{{\cal C}_a}$,
then
$y
\in 2Q_k \subset \Omega$ and 
$(y,t) \in {\cal E}_1$.
Indeed, since
$F$ is closed, there exists
$x
\in Q_k
\cap F$ such that
$|y-x|\le\frac t{20a}$.  Since
$t\le 10la$, we have
$|y-x|\le \frac l2$ and  the first claim is
true. Moreover, $|y-x| \le \frac t{20a} <t$, hence 
the second claim.

It follows in particular that  $\overline{{\cal C}_a}$ remains far
from the  boundary of $ \Omega \times ]0,\infty[$, so that we
do not care about the boundary values of $u_t(y)$, and that 
$|u_t(y)| \le \gamma\lambda$ on $\overline{{\cal C}_a}$.

 The Green-Riemann formula 
shows  
that
$I_{a}$ is equal to
\[
  -\int_{{\cal C}_{a}} 
\partial_tu_t(y)\,\overline{u_t(y)} \, dydt 
+\int_{\partial {\cal C}_{a}}
tB\overline\nabla u_t(y)\cdot
N_a(y,t) \,\overline{u_t(y)}\,d\sigma_{a}(y,t). 
\]
In this computation, $N_a(y,t)$ is the unit normal vector
outward
${\cal C}_{a}$  whereas $d\sigma_{a}$ is the surface
measure over $\partial {\cal C}_{a}$. 
Moreover, the Green-Riemann formula again yields
\[
2\mbox{ Re }\int_{{\cal C}_{a}}
\partial_tu_t(y)\overline{u_t(y)} dydt = \int_{\partial {\cal
C}_{a}} 
\left\vert u_t(y)\right\vert^2 N_a(y,t)\cdot (0,\ldots, 0,1)\,
d\sigma_{a}(y,t).
\]
Finally, 
\[
\begin{array}{lll}
\displaystyle \int_{{\cal C}_{a}} t\left\vert
\overline{\nabla}u_t(y)\right\vert^2 dydt  & \leq & 
 \displaystyle C \int_{\partial {\cal C}_{a}} 
\left\vert u_t(y)\right\vert^2 d\sigma_{a}(y,t) \\
\\
& + & \displaystyle C \int_{\partial {\cal C}_{a}}
t\left\vert u_t(y)\right\vert \left\vert \overline{\nabla}
u_t(y)\right\vert d\sigma_{a}(y,t).
\end{array}
\]

Since $\left\vert u_t(y)\right\vert \leq \gamma\lambda$ on 
$\partial {\cal C}_{a}$, we obtain that 
$$
\int_1^2 \int_{\partial {\cal C}_{a}} 
\left\vert u_t(y)\right\vert^2 d\sigma_{a}(y,t)da \le
\gamma^2\lambda^2 
\int_1^2 \int_{\partial {\cal C}_{a}} 
 d\sigma_{a}(y,t)da.
$$
We claim that
$$
\int_1^2 \int_{\partial {\cal C}_{a}} 
 d\sigma_{a}(y,t)da  \le c |C|. 
$$
Indeed, this integral is bounded by  $c\int_{\cal G} \frac{dzds}{s}$ 
where $\cal G$ is
the union of the sets $\partial {\cal C}_{a}$ for $1<a<2$.
This is the set of points $(z,s)$ with $z\in C$ and 
$\ep < s< 2\ep$ or $\psi(z) <s<2\psi(z)$ or $10l < s <20l$.
The claim follows readily.

It remains to establish
$$
\int_1^2\int_{\partial {\cal C}_{a}}
t\left\vert u_t(y)\right\vert \left\vert \overline{\nabla}
u_t(y)\right\vert d\sigma_{a}(y,t)da \le c\gamma^2\lambda^2|C|.
$$
Using the previous notation and a change of variables, this
integral is bounded by 
$$
\gamma\lambda \int_{\cal G} \left\vert \overline{\nabla}
u_t(y)\right\vert \,{dydt}. 
$$
Pick a  covering of ${\cal G}$ with bounded overlap by balls
$B_j=B((x_j,t_j), \frac{\ep t_j}{20})$. Remark  that
$(x,t)
\in B_j$  implies $t\sim t_j \sim r(B_j)$, the radius of $B_j$. 
 Then using  H\"older's
inequality and again Caccioppoli's inequality 
\[
\begin{array}{lll}
\displaystyle  \int_{\cal G} \left\vert
\overline{\nabla}  u_t(y)\right\vert \,{dydt} &\leq &\displaystyle 
c\sum
\int_{B_j} |\overline{\nabla}u_t(y)| \, dydt  
\\
\\
 &\le & \displaystyle  c\sum |B_j|^{1/2}
r(B_j)^{-1}
\left(\int_{2B_j} |u_t(y)|^2 \, dydt\right)^{1/2}
\\
\\
& \le& \displaystyle  c\gamma\lambda\sum |B_j|
r(B_j)^{-1}
\\
\\
& \le &  \displaystyle  c\gamma\lambda \int_{\widetilde{\cal  G}}
\frac{dzds}{s}
\\
\\
& \le & \displaystyle  c\gamma\lambda |C|.
\end{array}
\]
Here, $ \widetilde{\cal  G}$ is a set like $\cal  G$ but slightly
enlarged: it is contained set of points $(z,s)$ with $z\in C$ and 
$\ep/2 < s< 4\ep$ or $\psi(z)/2 <s<4\psi(z)$ or $5l < s <40l$. 
\rule{2mm}{2mm}
 
\begin{remark}
When $L$ is the Laplacian, then $(y,t) \mapsto u_t(y)$ is harmonic
so that the Caccioppoli inequality can be replaced  by the mean value
property and one can proceed directly using the square functions
(and not their variants).

When $\Omega\ne \reel^n$,  the Whitney cubes $Q_k$ are designed to
stay away from the boundary of $\Omega$ so that interior estimates
suffice. When  $\Omega = \reel^n$, one can proceed directly and
prove the good lambda inequality on $\reel^n$.   Details are left to
the reader. 
\end{remark}

\subsection{Carleson measure estimates and $BMO$}
\label{BMOCarleson}

In the present section, let $\Omega=\reel^n$ or $\Omega$ be a strongly
Lipschitz domain of $\reel^n$. Consider an operator $L=(A, \Omega,
V)$ satisfying (\ref{Gauss}) and the restrictions of Theorem
\ref{comparison} on $\Omega$.  Set
$P_{t}=e^{-tL^{1/2}}$. We intend to show that, when a function
$\phi$ belongs to
$BMO_{z,a}(\Omega)$ (resp. $BMO_{CW}(\Omega)$),
$\displaystyle t\left\vert \partial_{t}P_{t}\phi(x)\right\vert^{2}
{dxdt}$ is a Carleson measure when $V=W^{1,2}_{0}(\Omega)$
(resp.
$V=W^{1,2}(\Omega)$), with the obvious modifications when $\Omega=\reel^n$.
When $Q$ is a cube with center $x_Q \in \Omega$ and radius 
$r =\ell(Q)/2$, define the tent over $Q$ by
\[
T(Q)=\left\{(y,t)\in \Omega\times \left]0,+\infty\right[;\ \left\vert
y-x_{Q}\right\vert<r-t\right\}.
\]
We recall that $|\ |$ is the sup norm on $\reel^n$.
When $\phi$ is locally integrable on $\Omega$, set
\[
T\phi(x)=\left(\sup\limits_{Q\ni x} \frac 1{\left\vert Q\cap
\Omega\right\vert}
\int_{T(Q)} \left\vert \partial_{t}P_{t}\phi(y)\right\vert^{2}
t{dydt}\right)^{1/2},
\]
where the supremum is taken over all the cubes $Q$ centered in
$\Omega$ and containing $x$.

We will need a few facts from  functional calculus for the operator
$L=(A,\Omega,V)$ (see
\cite{mac} and
\cite{yagi}). 

First note that $L$ is
one-to-one (except if $\Omega$ is bounded and $V=W^{1,2}(\Omega)$,
which is excluded from our discussion), and if
$\omega=\sup\limits_{x\in\Omega,\ \xi\in {\bf C}^n} \left\vert \arg
A(x)\xi .\overline{\xi}\right\vert$, one has $\omega<\pi/2$   and $L$
is $\omega$-accretive on $V$ (see \cite{asterisque}). If $\mu\in
\left]\omega,\pi\right[$, and $\Gamma_{\mu}=\left\{z\in {\bf C}\setminus
\left\{0\right\}, \left\vert \arg z\right\vert <\mu\right\}$, for all
function $f\in H^{\infty}(\Gamma_{\mu})$, one can define a bounded operator
$f(L)$ on $L^2(\Omega)$ such that $\left\Vert f(L)\right\Vert \leq c_{\mu}
\left\Vert f\right\Vert_{H^{\infty}(\Gamma_{\mu})}$.

If $\psi\in H^{\infty}(\Gamma_{\mu})$ and if there exist $c,s>0$ such that,
for all $\zeta\in \Gamma_{\mu}$, $\left\vert \psi(\zeta)\right\vert \leq
c\left\vert \zeta\right\vert^s \left(1+\left\vert
\zeta\right\vert\right)^{-2s}$, then $\psi(L)$ may be computed thanks to
the Cauchy formula:
\[
\psi(L)=\frac 1{2i\pi} \int_{\gamma} \left(\zeta-L\right)^{-1}\psi(\zeta)d\zeta
\]
where $\gamma$ is made of two rays $re^{\pm i\nu}$, $r>0$, $\omega<\nu<\mu$
and is described counterclockwise.
If $\psi$ is such a function and is not identically zero, and if one
defines $\psi_t(\zeta)=\psi(t\zeta)$, there exists $c_{\psi}>0$ such that,
for all $f\in L^2(\Omega)$,
\begin{equation} \label{calculus}
c_{\psi}\left\Vert f\right\Vert_{2} \leq \left(\int_0^{+\infty}
\left\Vert \psi_t(L)f\right\Vert^2_2 \frac{dt}t\right)^{1/2} \leq
c_{\psi}^{-1} \left\Vert f\right\Vert_{2}.
\end{equation}
See \cite{mac}. This remark applies to $\psi(z)= -z^{1/2}
e^{-z^{1/2}}$ and $\psi_t(L)= t\partial_{t}P_{t}$.

The result of this section is a follows (see \cite{feffstein},
Theorem 3, $(iii)$, p. 145):
\begin{proposition} \label{Carleson}
Assume that $L=(A,\Omega,V)$ satisfies (\ref{Gauss}).
\begin{itemize}
\item[$(a)$]
If $\Omega=\reel^n$, then there exists $C>0$ such that, for all
$\phi\in BMO(\reel^n)$, $\left\Vert T\phi\right\Vert_{L^{\infty}(\reel^n)}
\leq C\left\Vert \phi\right\Vert_{BMO(\reel^n)}$.
\item[$(b)$]
If  $L$ satisfies DBC, then there exists
$C>0$ such that, for all
$\phi\in BMO_{z,a}(\Omega)$, $\left\Vert T\phi\right\Vert_{L^{\infty}(\Omega)}
\leq C\left\Vert \phi\right\Vert_{BMO_{z,a}(\Omega)}$.
\item[$(c)$]
If $\Omega$ is unbounded and $L$ satisfies NBC, then there exists
$C>0$ such that, for all
$\phi\in BMO_{CW}(\Omega)$,
$\left\Vert T\phi\right\Vert_{L^{\infty}(\Omega)} \leq C\left\Vert
\phi\right\Vert_{BMO_{CW}(\Omega)}$.
\end{itemize}
\end{proposition}
The proofs of the three assertions are similar, and we give the one of
$(b)$. Consider $\phi\in BMO_{z,a}(\Omega)$ with $\left\Vert
\phi\right\Vert_{BMO_{z,a}(\Omega)}\leq 1$, 
$x\in \Omega$ and a cube $Q$
centered
in $\Omega$ and containing $x$.
 Write
\[
\phi=\phi_{2Q \cap \Omega} + (\phi-\phi_{2Q \cap \Omega}).
\]
It is classical using the square function estimate for
$\partial_{t}P_{t}$,  the decay of the kernel of
$\partial_{t}P_{t}$ and $BMO$ inequalities that 
\[
\frac 1{\left\vert Q\cap
\Omega\right\vert}
\int_{T(Q)} \left\vert \partial_{t}P_{t}(\phi-\phi_{2Q \cap
\Omega})(y)\right\vert^{2} t{dydt} \le c
\|\phi\|_{BMO_{CW}(\Omega)}^2
\]
and since $BMO_{z,a}(\Omega) \subset BMO_{r,a}(\Omega)=
BMO_{CW}(\Omega)$, $\|\phi\|_{BMO_{CW}(\Omega)} \le c$. 

It remains to control 
$I_Q=\frac 1{\left\vert Q\cap
\Omega\right\vert}
\int_{T(Q)} \left\vert \partial_{t}P_{t}(\phi_{2Q \cap
\Omega})(y)\right\vert^{2} t{dydt}$ by a constant which does not
depend on $Q$.

We use the following lemma.

\begin{lemma} \label{mean}
Let $\phi\in BMO_{z,a}(\Omega)$ with $\left\Vert
\phi\right\Vert_{BMO_{z,a}(\Omega)}\leq 1$ and $Q$ be a cube centered
in $\Omega$. Then 
$$\left\vert
\phi_{Q
\cap
\Omega}\right\vert
\leq C \sup (\ln\left(\frac{\delta}r\right), 1),$$ where  $r$ is
the radius of $Q$ and $\delta\ge 0$ its distance to the boundary of
$\Omega$.
\end{lemma}

{\bf Proof:} Notice first that, for any cube $Q$ of type $(a)$, one has
\[
\left\vert \phi_{Q}-\phi_{2Q}\right\vert \leq C.
\]
Moreover, if $Q$ is of type $(b)$, then $|\phi_Q| \le 1$ by definition
of
$BMO_{z,a}(\Omega)$ and H\"older's inequality.

 Assume first that $Q$ is of type $(a)$, and let $k$
be the smallest integer such that $2^{k}Q$ is of type $(b)$. Then, one
has
\[
\left\vert\phi_Q\right\vert  \leq  
\sum\limits_{i=0}^{k-1} \left\vert
\phi_{2^{i}Q}-\phi_{2^{i+1}Q}\right\vert + \left\vert
\phi_{2^kQ}\right\vert  \leq  C(k+1)  \leq  C^{\prime}
\ln\left(\frac{\delta}r\right).
\]
Assume now that $4Q\cap\partial\Omega\neq \emptyset$. Take a  Whitney
decomposition of $Q\cap\Omega$ with respect to $\partial\Omega$,
\[
Q\cap\Omega=\bigcup\limits_k Q_k,
\]
where, for each $k$, $Q_k$ is a type $(b)$ cube. Therefore, one has
\[
 \left\vert \phi_{Q\cap\Omega}\right\vert  \leq 
 \sum\limits_k \frac{\left\vert Q_k\right\vert}{\left\vert
Q\cap\Omega\right\vert} \left\vert \phi_{Q_k}\right\vert 
 \leq  1.
\]

Let us come back to the proof of Proposition \ref{Carleson}.  
Lemma \ref{poissint} (see Appendix A) shows
that, for all
$y\in Q \cap \Omega$,
\[
 \left\vert \int_{\Omega} \partial_{t}
p_t(z,y)dz\right\vert 
\leq   \frac Ct
\left(1+\frac{ \delta(y)}{t}\right)^{-1}
\]
where $\delta(y)$ is the distance from $y$ to the boundary  of
$\Omega$. It is then fairly easy using that $\Omega$ is strongly
Lipschitz  to show that
\[
\frac 1{\left\vert Q\cap
\Omega\right\vert}
\int_{T(Q)}  \frac C{t^2}
\left(1+\frac{ \delta(y)}{t}\right)^{-2}
t{dydt}
 \leq C \inf (\frac {r^2}{\delta^2}, 1).
\]
See \cite{autchdom}, Lemma 29, for details in a related situation.
This and Lemma \ref{mean} prove that $I_Q \le C$.
Assertion $(b)$ is complete.

For assertions $(a)$ and $(c)$, decompose $\phi$ as
above and since $\partial_{t}P_{t}$ annihilates constants, only the
first term arises. Proposition
\ref{Carleson} is proved. \rule{2mm}{2mm}

\subsection{A weakly dense class}

We shall need the
following lemma (see \cite{grenoble}, Lemme II.11):

\begin{lemma} \label{density} Assume that $L=(A,\Omega,V)$ satisfies
$(G_{\infty})$.
For all function $f\in H^1_{max,L}(\Omega)$, there exists a sequence
$(f_k)_{k\in
\nat}$ of functions in
$H^1_{max,L}(\Omega)\cap L^2(\Omega)$, such that, for all $\phi\in
C_c(\Omega)$,
\[
\lim\limits_{k\rightarrow +\infty} \langle f_k,\phi\rangle = \langle
f,\phi\rangle
\]
and, for all $k\in \nat$,
\[
\left\Vert f_{k} \right\Vert_{H^{1}_{max,L}(\Omega)} \leq \left\Vert
f\right\Vert_{H^{1}_{max,L}(\Omega)}.
\]
\end{lemma}

We  briefly sketch the proof. Let $f\in
H^{1}_{max,L}(\Omega)$. Define
 $f_k=P_{2^{-k}}f$. 
Then, the decay of the kernel of $P_{2^{-k}}$ and $f \in
L^1(\Omega)$ imply
$f_k\in L^2(\Omega)$. It is obvious that 
$\left\Vert f_k\right\Vert_{H^1_{max,L}(\Omega)} \leq \left\Vert
f\right\Vert_{H^1_{max,L}(\Omega)}$.
Using arguments analogous to \cite{grenoble}, p. 776,  which rely
on the decay of the kernel of $P_{2^{-k}}$, one obtains the weak
convergence.
\rule{2mm}{2mm}

\subsection{Proof of
Theorem \ref{comparison}}
\label{comparglobal}
We begin this section with the proof of assertion $(a)$ in Theorem
\ref{comparison}.

\paragraph{$H^{1}(\reel^n)\subset
H^{1}_{max,L}(\reel^n)$:}
If $a$ is an $H^{1}(\reel^n)$-atom, it is plain to see, using
$(G_{\infty})$, that $\left\Vert a^{\ast}_L\right\Vert_{1}\leq C$ (see,
for instance, \cite{garcia}).

\paragraph{$H^{1}_{max,L}(\reel^n)\subset
H^{1}(\reel^n)$:}
Consider first $f\in H^1_{max,L}(\reel^n)\cap L^2(\reel^n)$. Observe
that, for all $z$ in
a sector $\Gamma_{\mu}$ with $\mu\in \left]\omega,\pi\right[$,
\[
\int_0^{+\infty}
\left(-tz^{1/2}e^{-tz^{1/2}}\right)\left(-tz^{1/2}e^{-tz^{1/2}}\right)\frac{dt}t
=\frac 14.
\]
As a consequence, one has
\[
\mbox{Id}=4\int_0^{+\infty} (tL^{1/2} P_t)(tL^{1/2} P_t)\frac{dt}t,
\]
where the integral converges strongly in $L^2(\reel^n)$. Note that 
$tL^{1/2} P_t= - t\partial_tP_t$. Thus,
if 
$f \in L^2(\reel^n)$ and $\phi$ is continuous with compact support in
$\reel^n$,  
 one has
\begin{equation} \label{functcalc}
\int_{\reel^n} f(y)\overline{\phi(y)}dy=4\int_{\reel^n}\int_0^{+\infty}
 (t\partial_t P_t)(f)(y) \overline{\left(t\partial_t P_t^{\ast}\right)
(\phi)(y)}\ \frac{dydt}t.
\end{equation}
We want to show that
\begin{equation} \label{inequality}
\left\vert \int_{\reel^n} f(y)\overline{\phi(y)}dy\right\vert \leq
C\left\Vert f\right\Vert_{H^1_{max,L}(\reel^n)} \left\Vert
\phi\right\Vert_{BMO(\reel^n)}.
\end{equation}
Well-known arguments from the theory of tent spaces 
(see \cite{tent, real}) show that the right-hand side of 
(\ref{functcalc}) is bounded by
\[
C\left\Vert sf\right\Vert_{L^1(\reel^n)} \left\Vert
T\phi\right\Vert_{L^{\infty}(\reel^n)},
\]
where  $T\phi$ is defined in Section
\ref{BMOCarleson} with $L$ replaced by $L^*$ and 
$$
sf(x) = \left(\int_{\Gamma_{1}(x)}
t^{1-n}\left\vert
\partial_tP_tf(y)\right\vert^2 dydt\right)^{1/2}
.$$ 
On the one hand, it is clear that $sf\le S_1f$ where $S_1$ is defined
in Section
\ref{area} and
Proposition
\ref{areaintegral} shows that
\[
\left\Vert S_1f \right\Vert_{L^1(\reel^n)}\leq C\left\Vert
f\right\Vert_{H^{1}_{max,L}(\reel^n)}.
\]
On the other hand,  assertion $(a)$ in Proposition \ref{Carleson} yields
\[
\left\Vert T\phi\right\Vert_{L^{\infty}(\reel^n)} \leq C\left\Vert
\phi\right\Vert_{BMO(\reel^n)}.
\]
This ends the proof of (\ref{inequality}). 

Up to now, we have proved that, when $f\in L^2(\reel^n)\cap
H^1_{max,L}(\reel^n)$,
\[
\left\Vert f\right\Vert_{H^1(\reel^n)} \leq
C\left\Vert f\right\Vert_{H^1_{max,L}(\reel^n)}.
\]
 
Consider $f\in
H^1_{max,L}(\reel^n)$.
Let $(f_k)_{k\in \nat}$ be a sequence given by
Lemma \ref{density}. Then, for all $k\in \nat$, since $f_{k}\in
H^{1}_{max,L}(\reel^n)\cap L^{2}(\reel^n)$, one has
\[
\left\Vert f_k\right\Vert_{H^1(\reel^n)} \leq C \left\Vert
f_k\right\Vert_{H^1_{max,L}(\reel^n)} \leq C\left\Vert
f\right\Vert_{H^1_{max,L}(\reel^n)}.
\]
Since $H^1(\reel^n)$ is the dual of $VMO(\reel^n)$, there exists a
subsequence
$(f_{\phi(k)})_{k\in \nat}$ and a function $g\in H^1(\reel^n)$
such that, for all $\phi\in C_c(\reel^n)$,
$\langle f_{\phi(k)}, \phi\rangle \rightarrow \langle g,\phi\rangle$. Since
Lemma \ref{density} shows that
$\langle f_{\phi(k)}, \phi\rangle \rightarrow \langle f,\phi\rangle$, one
has $f=g$. Therefore, $f\in H^1(\reel^n)$ and
\[
\left\Vert f\right\Vert_{H^1(\reel^n)} \leq \underline{\lim} \left\Vert
f_{\phi(k)}\right\Vert_{H^1(\reel^n)} \leq C\left\Vert
f\right\Vert_{H^{1}_{max,L}(\reel^n)}.
\]

\bigskip

We now turn to the proof of assertion $(b)$ in Theorem \ref{comparison}.

\paragraph{$H^{1}_{r,a}(\Omega)\subset
H^{1}_{max,L}(\Omega)$:}
Let $a$ be an atom of type $(a)$ supported in a cube $Q$.
Then, using $(G_{\infty})$, one sees that $\left\Vert
a\right\Vert_{H^1_{max,L}(\Omega)} \leq C$ (see \cite{garcia}). If $a$ is
of type $(b)$, write that
$P_ta(x)=\int_Q \left[p_t(x,y)-p_t(x,y_0)\right]a(y)dy$ where
$y_0$ is a point on $\partial \Omega$ such that $\left\vert
y-y_0\right\vert \sim d(y,\partial \Omega)$ whenever $y \in \supp a$
(remember that
$p_{t}(x,y_0)=0$ since $y_0\in \partial \Omega)$ and use this
representation and
$(G_{\infty})$ to
show that $\left\Vert a\right\Vert_{H^1_{max,L}(\Omega)} \leq C$.

\paragraph{$H^{1}_{max,L}(\Omega)\subset
H^{1}_{r,a}(\Omega)$:}
As in the case of $\reel^n$, consider first $f\in H^1_{max,L}(\Omega)\cap
L^2(\Omega)$.
Arguing as before, one obtains that, if $\phi$ is continuous with compact
support, one has
\begin{equation} \label{functcalcdomain}
\int_{\Omega} f(y)\overline{\phi(y)}dy=4\int_{\Omega}\int_0^{+\infty}
 (t\partial_t P_t)(f)(y) \overline{\left(t\partial_t P_t^{\ast}\right)
(\phi)(y)}\ \frac{dydt}t.
\end{equation}
We want to show that
\begin{equation} \label{inequalitydomain}
\left\vert \int_{\Omega} f(y)\overline{\phi(y)}dy\right\vert \leq
C\left\Vert f\right\Vert_{H^1_{max,L}(\Omega)} \left\Vert
\phi\right\Vert_{BMO_{z,a}(\Omega)}.
\end{equation}
Use the theory of tent spaces again (see \cite{tent}) to obtain
\[
\left\vert \int_{\Omega} f(y)\overline{\phi(y)}dy\right\vert \leq
C\left\Vert S_1f\right\Vert_{L^1(\Omega)} \left\Vert
T\phi\right\Vert_{L^{\infty}(\Omega)}.
\]
Proposition \ref{areaintegral} shows that
\[
\left\Vert S_1f \right\Vert_{L^1(\Omega)}\leq C\left\Vert
f\right\Vert_{H^{1}_{max,L}(\Omega)},
\]
whereas assertion $(b)$ in Proposition \ref{Carleson} yields
\[
\left\Vert T\phi\right\Vert_{L^{\infty}(\Omega)} \leq C\left\Vert
\phi\right\Vert_{BMO_{z,a}(\Omega)},
\]
which ends the proof of (\ref{inequalitydomain}). This inequality implies
that there exists $C>0$ such that, for all $f\in H^1_{max,L}(\Omega)\cap
L^2(\Omega)$,
\begin{equation} \label{partialresult}
\left\Vert f\right\Vert_{H^1_{r,a}(\Omega)} \leq C\left\Vert
f\right\Vert_{H^1_{max,L}(\Omega)}.
\end{equation}
Indeed,  since $BMO_{z}(\Omega)=BMO_{z,a}(\Omega)$ and  
$H^1_{r,a}(\Omega)=H^1_{r}(\Omega)= (VMO_{z}(\Omega))'$, we have that
$\left\Vert f\right\Vert_{H^1_{r,a}(\Omega)}\sim \sup\{
|\langle f, \phi \rangle|; \phi \in C_c(\Omega),
\left\Vert\phi\right\Vert_{BMO_{z,a}(\Omega)}=1\}$
 (see
Section \ref{BMO}).

One gets rid of the
condition $f\in L^{2}(\Omega)$ using
Lemma
\ref{density}. Indeed, if $f\in
H^1_{max,L}(\Omega)$, there exists a
sequence $(f_k)_{k\in
\nat}$ of functions in
$H^1_{max,L}(\Omega)\cap L^2(\Omega)$, such that, for all $\phi\in
C_c(\Omega)$,
\begin{equation} \label{prop1}
\lim\limits_{k\rightarrow +\infty} \langle f_k,\phi\rangle = \langle
f,\phi\rangle
\end{equation}
and, for all $k\in \nat$,
\begin{equation} \label{prop2}
\left\Vert f_{k} \right\Vert_{H^{1}_{max,L}(\Omega)} \leq C\left\Vert
f\right\Vert_{H^{1}_{max,L}(\Omega)}.
\end{equation}
By (\ref{partialresult}) and (\ref{prop2}) and $H^1_{r}(\Omega)=
H^1_{r,a}(\Omega)$, the
$f_k$'s are bounded in
$H^1_{r}(\Omega)$. Since $H^1_{r}(\Omega)$ is the dual of
$VMO_z(\Omega)$, there exists a subsequence
$(f_{\phi(k)})$ which converges $*$-weakly to $g\in H^1_{r}(\Omega)$.
Then, (\ref{prop1}) implies that $f=g$. Moreover,
\[
\left\Vert f\right\Vert_{H^1_{r,a}(\Omega)} \sim \left\Vert
f\right\Vert_{H^1_{r}(\Omega)} \leq C\left\Vert
f\right\Vert_{H^{1}_{max,L}(\Omega)}.
\]

\bigskip

We are now left with the task of proving assertion $(c)$ in Theorem
\ref{comparison}. We shall  prove that $H^{1}_{z}(\Omega)\subset
H^1_{max,L}(\Omega)\subset H^{1}_{CW}(\Omega)$ since
$H^{1}_{CW}(\Omega)= H^{1}_{z,a}(\Omega)\subset H^1_z(\Omega)$ by
Theorem \ref{H1equalities}.

\paragraph{$H^{1}_{z}(\Omega)\subset
H^1_{max,L}(\Omega)$:}
Recall that $p_t(x,y)$ is the Poisson kernel for $L$. Recall that 
(G$_\infty$) and the subordination formula imply that 
$$
|p_t(x,y)| \le C t^{-n} \left( 1+
\frac{|x-y|}{t}\right)^{-n-1}
$$
and
$$|p_t(x,y) -p_t(x,y')| \le  C t^{-n}\left(\frac{\left\vert
y-y^{\prime}\right\vert}{t}\right)^{\nu}
\left(\left( 1+
\frac{|x-y|}{t}\right)^{-n-1-\nu} +   \left( 1+
\frac{|x-y'|}{t}\right)^{-n-1-\nu}\right)
$$
for some $\nu \in (0, 1]$.
For all $t>0$ and $x\in \Omega$, define
\[
F_{x,t}(y)=t^{n}\left( 1+
\frac{|x-y|}{t}\right)^{n+1}p_t(x,y).
\]
It is easy to show that 
\[
\left\vert F_{x,t}(y)\right\vert \leq C 
\]
and 
\[
\left\vert F_{x,t}(y)- F_{x,t}(y')\right\vert \leq C \left(\frac{\left\vert
y-y^{\prime}\right\vert}{t}\right)^{\nu} 
\]
for all $y, y'\in \Omega$.
Thus, the function $F_{x,t}$ may be extended to a bounded
H\"older continuous function on $\overline{\Omega}$, then on
$\reel^n$ (see
\cite{singular}, Chapter 6, p. 174, Theorem 3). If this extension is
denoted by $\widetilde{F}_{x,t}$, one has
\[
\left\vert \widetilde{F}_{x,t}(y)\right\vert \leq C_0C
\]
and
\[
\left\vert
\widetilde{F}_{x,t}(y)-\widetilde{F}_{x,t}(y^{\prime})\right\vert \leq
C_0C\left(\frac{\left\vert
y-y^{\prime}\right\vert}{t}\right)^{\nu},
\]
for all $y,y^{\prime}\in \reel^n$, where $C_0$ only depends on $\Omega$.
Define now
\[
\widetilde{p}_t(x,y)=t^{-n}\left( 1+
\frac{|x-y|}{t}\right)^{-n-1}\widetilde F_{x,t}(y).
\]
Then, one has
\begin{equation} \label{estim1}
\left\vert \widetilde{p}_t(x,y)\right\vert \leq 
C{t^{-n}} \left(1+
\frac{|x-y|}{t}\right)^{-n-1}
\end{equation}
and
\begin{equation} \label{estim2}
\left\vert
\widetilde{p}_t(x,y)-\widetilde{p}_t(x,y^{\prime})\right\vert
\leq  C{t^{-n}} \left(\frac{\left\vert
y-y^{\prime}\right\vert}{t}\right)^{\nu}
\end{equation}
for all $x\in \Omega$, $ y,y^{\prime} \in \reel^n$ and all
$t>0$. Moreover, for all
$t>0$ and all $x,y\in \Omega$,
$p_t(x,y)=\widetilde{p}_t(x,y)$.

Consider now a function $f\in H^1_z(\Omega)$, extended by $0$ outside
$\Omega$, so that $f\in H^1(\reel^n)$ and
$
\left\Vert f\right\Vert_{H^1_z(\Omega)}= \left\Vert
f\right\Vert_{H^1(\reel^n)}.
$
 For all $x\in \Omega$, one has
\begin{equation} \label{Ktilde}
\begin{array}{lll}
\displaystyle \int_{\Omega} p_t(x,y)f(y)dy = \displaystyle
\int_{\reel^n}\widetilde{p}_t(x,y)f(y)dy.
\end{array}
\end{equation}
Using the atomic decomposition of $f$ into $H^1(\reel^n)$-atoms and the
estimates (\ref{estim1}) and (\ref{estim2}) for $\widetilde{p_t}$, one
easily deduces from (\ref{Ktilde}) that  
\[
\left\Vert f^{\ast}_{L}\right\Vert_{L^1(\Omega)} \leq
C\left\Vert f\right\Vert_{H^1(\reel^n)}.
\] 

This ends the proof of the inclusion $H^1_z(\Omega)\subset
H^1_{max,L}(\Omega)$.

\paragraph{$H^{1}_{max,L}(\Omega)\subset
H^{1}_{CW}(\Omega)$:}
Arguing as in the proof of $H^{1}_{max,L}(\Omega)\subset
H^{1}_{r,a}(\Omega)$ under
DBC, one proves that, for all $f\in H^{1}_{max,L}(\Omega)\cap
L^{2}(\Omega)$ and all function
$\phi$ continuous and compactly supported in $\Omega$,
\[
\left\vert \int_{\Omega} f(x)\overline{\phi(x)}dx \right\vert \leq
C\left\Vert f\right\Vert_{H^{1}_{max,L}(\Omega)} \left\Vert
\phi\right\Vert_{BMO_{CW}(\Omega)}.
\]
The proof uses the theory of tent spaces, Proposition \ref{areaintegral}
and assertion $(c)$ of Proposition \ref{Carleson}.

Then, one gets rid of the assumption $f\in L^{2}(\Omega)$ as before,
using Lemma \ref{density} and the fact that 
$H^1_{CW}(\Omega)$ is the dual space of $VMO_{CW}(\Omega)$.

\subsection{Some consequences} 

Assume the hypotheses of Theorem \ref{comparison}.  We list some
consequences of the proofs.

Each maximal Hardy space is characterized by the square functions
$S_1f$ and $sf$ being in $L^1$.  Indeed, we have already seen that
$\|sf\|_1 \le \| S_1f\|_1 \le C\| f^*_L\|_1$ and the argument
via tent spaces and Carleson measures of Section 2.7 shows in
fact that
$sf
\in L^1$ implies
$f$ is in an atomic space. As the atomic space is contained in a
maximal space, we have a full circle of implications.  See the
forthcoming paper \cite{adm} where a general theory for Hardy
spaces defined through square functions of type $sf$ associated
to abstract operators $L$ is developed.

Each of the atomic $BMO$ space, has a characterization in terms of 
Carleson measures. In other words, the converse to the inequalities 
of Proposition
\ref{Carleson} between $BMO$ norms and $\|T\phi\|_\infty$ hold
(provided $\phi$ satisfies an integrability condition as in
\cite{feffstein}). We leave to the reader the care of checking this.

\subsection{``Easy'' embeddings between Hardy spaces}

We prove part of  Theorem \ref{H1equalities} for the global Hardy
spaces.

\paragraph {Assertion $(a1)$:} It is proved in
\cite{stkrch} (p. 305, proof of Theorem 2.7, $(1) \Rightarrow (2)$) for
local spaces and when $\Omega$ is bounded. We briefly
recall the argument for completeness. Let $f\in
H^1_r(\Omega)$ and $F\in H^1(\reel^n)$ be an extension of $f$
satisfying $\left\Vert F\right\Vert_{H^1(\reel^n)}\leq
2\left\Vert f\right\Vert_{H^1_r(\Omega)}$. The function $F$
may be decomposed into
\[
F=\sum\limits_Q \lambda_Qa_Q.
\]
In this sum, we are only interested in the cubes $Q$ which intersect
$\Omega$. If $4Q\subset \Omega$, consider $a_Q$ as a type $(a)$ atom. If
$2Q\subset \Omega$ and $4Q\cap \partial\Omega\neq \emptyset$, consider
$a_Q$ as a type $(b)$ atom. Finally, consider the case when $2Q\cap
\partial \Omega\neq \emptyset$ and perform a Whitney decomposition of
$Q\cap\Omega$ with respect to $\partial\Omega$:
\[
Q\cap \Omega=\bigcup\limits_k Q_k
\]
where each $Q_k$ is a type $(b)$ cube and decompose
\[
a_Q{\bf 1}_{\Omega}=\sum\limits_k \left\vert Q_k\right\vert^{1/2}\left\Vert
a_Q{\bf 1}_{Q_k}\right\Vert_2\frac{a_Q{\bf 1}_{Q_k}}{\left\vert
Q_k\right\vert^{1/2}\left\Vert a_Q{\bf 1}_{Q_k}\right\Vert_2}.
\]
Since each $Q_k$ is a type $(b)$ cube, $\frac{a_Q{\bf
1}_{Q_k}}{\left\vert Q_k\right\vert^{1/2}\left\Vert a_Q{\bf
1}_{Q_k}\right\Vert_2}$ is a type
$(b)$ atom. Moreover,
\[
\begin{array}{lll}
\displaystyle \sum\limits_k \left\vert Q_k\right\vert^{1/2}\left\Vert
a_Q{\bf 1}_{Q_k}\right\Vert_2& \leq & \displaystyle \left(\sum\limits_k
\left\vert Q_k\right\vert\right)^{1/2} \left(\sum\limits_k \left\Vert
a_Q{\bf 1}_{Q_k}\right\Vert_2^2\right)^{1/2} \\
\\
& \leq & \displaystyle \left\vert Q\right\vert^{1/2} \left\Vert
a_Q\right\Vert_2 \\
\\
& \leq & 1.
\end{array}
\]
Thus, we have obtained a decomposition of $f=\left. F\right\vert_{\Omega}$
into $H^1_{r,a}(\Omega)$-atoms.

\paragraph{Assertion $(a2)$:}  Let $a$ be an $H^1_{r,a}(\Omega)$-atom.
If
$a$ is of type $(a)$, then its extension by $0$ outside $\Omega$ is an
$H^1(\reel^n)$-atom. Hence $a \in H^1_r(\Omega)$. 

If $a$ is of type $(b)$, let $Q$ be a type $(b)$ cube on which $a$ is
supported. We use the following claim, whose proof is deferred to
Appendix B:

\noindent{\bf Claim:} There exists $\rho \in ]0,+\infty]$, such that 
if
$Q$ is a type
$(b)$ cube and $\ell(Q)< \rho$,   there exists a cube
$\widetilde{Q}\subset
{}^c\Omega$ such that 
$\left\vert \widetilde{Q}\right\vert \sim \left\vert Q\right\vert$ and 
the distance from $\widetilde{Q}$ to $Q$ is comparable 
to  $\ell(Q)$. Furthermore, $\rho=\infty$ is $^c\Omega$ is
unbounded.

Define an extension of $a$ as follows: Let $A(x) = a(x) $ if $x\in Q$,
$A(x)= -\frac 1{|Q|} \int_Q a$ if $x \in \tilde Q$. Let $Q_0$ be the
smallest cube in $\reel^n$ containing $Q$ and $\tilde Q$.  It is clear
that 
$\supp A \subset Q_0 $, $\|A\|_2 \le C|Q_0|^{-1/2}$,  that $\int
A=0$ and that $a$ is the restriction of $A$ to $\Omega$. Hence, $a \in
H^1_r(\Omega)$ 

\paragraph{Assertion $(b1)$:} The inclusion $H^1_{z,a}(\Omega) \subset
H^1_{CW}(\Omega)$ is obvious. We give the converse using an argument due to 
Lou and McIntosh.

Let $a$ be an $H^1_{CW}(\Omega)$-atom associated to a cube $Q$. We want to show
that
$a$ belongs to $H^1_{z,a}(\Omega)$. In fact, we are going to show that $a$ can
be written as a sum of type $(a)$ atoms.

If $a$ is supported on a type $(a)$ cube, we do nothing. If not, $a$ is
supported by $Q \cap \overline\Omega$ where $Q$ is a cube centered in
$\Omega$ which is not of type $(a)$. Since $a$ is square integrable with
mean value 0 on the Lipschitz domain
$Q\cap \Omega$, we can invoke the following corollary of a  result by Ne\v cas
\cite{necas}, Chapter 3, Lemma 7.1.

\begin{lemma} Let $D$ be a bounded Lipschitz domain. The divergence operator
is  a (continuous) map from $H^1_0(D)^n$ onto $L^2_0(D)=\{f\in L^2(D); \int_D
f = 0\}$: there exists $C>0$ depending only on the Lipschitz constant of $D$
such that for all $f \in L^2_0(D)$, there exists ${\bf g} \in H^1_0(D)^n $
such that $div {\bf g} =f$ and $ \int_D 
|\nabla {\bf g}|^2 \le C \int_D
 |f|^2$.
\end{lemma}  

Indeed,  Ne\v cas proves that the gradient operator is one-one
with closed range from $L^2_0(D)$ into $H^{-1}(D)^n$ with $\|\nabla
f\|_{H^{-1}(D)^n} \ge C \|f\|_2$ and controls the constant.  

Hence, pick 
 $b\in H^1_0(Q\cap \Omega)^n$ with $a= div b$ and $ \int 
|\nabla b|^2 \le C \int |a|^2$, the constant
$C$ depending only on the Lipschitz constant of $Q\cap
\Omega$, henceforth only on $\Omega$. Extend
$b$ by 0 outside of $Q\cap \Omega$. 

Pick a Whitney
decomposition $(Q_k)$ of $\Omega$
by cubes 
 from a dyadic grid  (of $\reel^n$) containing
$Q$ and so that  $8Q_k \subset \Omega$. Again
\[
Q\cap \Omega=\bigcup\limits_{k\in K} Q_{k}
\]
and 
\[
\sum\limits_{k \in K} |Q_k| =|Q\cap \Omega|.
\]
for some index set $K$.
Let $(\eta_k)$ be a smooth partition of unity associated to this covering
with $\eta_k$ supported in $2Q_k$ and $\|\eta_k\|_\infty\le 1$ and 
$\|\nabla \eta_k\|_\infty \le C(\ell(Q_k))^{-1}$.
Write $a =\sum_{k\in K} div (\eta_k b)$ and set 
$$a_k = \frac{div (\eta_k
b)}{\|div (\eta_k b)\|_2 |2Q_k|^{1/2}} \quad  {\rm and} \quad \lambda_k=
{\|div (\eta_k b)\|_2 |2Q_k|^{1/2}}$$ whenever this number is not zero.
Otherwise set
$\lambda_k=0$ and $a_k=0$. It is clear from its construction that 
$a_k$ is a type $(a)$ atom. It remains to show that $\sum \lambda_k \le
C$
 independent of $Q$.  By Cauchy-Schwarz inequality this sum does not
exceed 
$$\left( \sum_{k\in K} 2^n|Q_k|\right)^{1/2} 
\left( \sum_{k\in K} \int_{2Q_k} |div (\eta_k b)|^2 \right)^{1/2}
$$
so it suffices to establish that the second term $\Sigma$ is controlled by
$\|a\|_2+
\|\nabla b\|_2$.

To this end, write  
$$
\Sigma \le \left( \sum_{k\in K} \int_{2Q_k} |\eta_k div b|^2
\right)^{1/2} + \left( \sum_{k\in K} \int_{2Q_k} |\nabla \eta_k\cdot b|^2
\right)^{1/2}.
$$
 For the term containing $div b = a$ using the finite covering property
of the cubes $2Q_k$ which leads to the bound $C\|a\|_2$. For 
the other term, observe that when $x \in 2Q_k$ then $d(x,\partial\Omega)
\sim \ell(Q_k)$. Hence, this and the finite overlap property of the cubes
$2Q_k$ lead to
$$\left( \sum_{k\in K} \int_{2Q_k} |\nabla \eta_k\cdot b|^2
\right)^{1/2}
 \le C \left( \sum_{k\in K} \int_{2Q_k}\left| \frac{b(x)
}{d(x,\partial\Omega)} \right|^2 dx
\right)^{1/2} \le  C\left(\int_{\Omega}\left| \frac{b(x)
}{d(x,\partial\Omega)}\right|^2 dx\right)^{1/2}.
$$
Now use  Hardy's inequality (see, e.g., \cite{davies}, Chapter 1, Section 5)
$$\left(\int_{\Omega}\left| \frac{b(x)
}{d(x,\partial\Omega)}\right|^2 dx\right)^{1/2} \le
\left(\int_{Q\cap\Omega}\left| \frac{b(x) }{d(x,\partial(Q\cap
\Omega))}\right|^2 dx\right)^{1/2} \le C
\|\nabla b \|_2$$
since $b \in H^1_0(Q\cap \Omega)$ and $Q\cap \Omega$ is strongly Lipschitz
and bounded. Note that the constant
$C$ in this inequality  depends only on the domain $\Omega$ and not, in
particular,  on the size of $Q$ (by a scaling argument). This ends the proof
of assertion
$(b1)$.

\paragraph{Assertion $(b2)$:}  It is obvious that 
$H^1_{z,a}(\Omega) \subset  H^1_{z}(\Omega) $.

\begin{remark}Let us see that Theorem \ref{H1equalities}, (a2),
is sharp.  Let us take   $\Omega=\reel
\setminus [0,1]$ (although $\Omega$ is not connected: the calculations
are the same in any dimension). Let
$a=-(2N)^{-1} {\cal X}_{[N,3N]}$ where $N$ is a positive integer.
Then $a$ is a type $(b)$ atom.  Suppose that the continuous embedding
$H^1_{r,a}(\Omega) \subset H^1_{r}(\Omega)$ holds. Then
$\|a\|_{H^1_r(\Omega)} \le C$ uniformly with respect to
$N$. Pick
$A\in H^1(\reel)$ with norm not exceeding $C+1$  and whose
restriction to
$\Omega$ is $a$. Since $\int_{\reel} A=0$ it must be that
$\int_{[0,1]} A=1$. 
Observe that 
$$
\|A\|_{H^1(\reel)} \ge c\int_2^{N/4} \sup_{2<t<N/2} |\phi_t*A(x)| \, dx
$$
by the maximal definition of $H^1(\reel)$ with $\phi$ smooth, supported
in $[0,1]$ and $\phi(x)>0$ for $x \in ]0,1[$.
Write for $x \in [2, N/4]$, and $t\in [2, N/2]$
$$\phi_t*A(x)= \int_0^1 (\phi_t (x-y) -\phi_t(x)) A(y) \, dy +
\phi_t(x). $$
Then,  for all $x$ in that range, the support condition of $\phi$
implies that  $\sup_{2<t<N/2} |\phi_t(x)| \ge C/x$
while   
$\sup_{2<t<N/2}|\int_0^1 (\phi_t (x-y) -\phi_t(x)) A(y) \, dy| \le
C/x^2$. We obtain therefore, 
$\|A\|_{H^1(\reel)} \ge C \ln N$, which is a contradiction.

\end{remark}

\subsection{Duality results and $BMO$ embeddings}\label{duality}

 We prove the parts of theorem \ref{thdual} left aside in Section
\ref{BMO}. 

\paragraph{Assertion $(a)$:} It is plain to see that,
if $\phi\in BMO_{z,a}(\Omega)$, $\phi$ defines a bounded linear
functional ${\cal L}$ on $H^{1}_{r,a}(\Omega)$ with $\left\Vert {\cal
L}\right\Vert
\leq \left\Vert \phi\right\Vert_{BMO_{z,a}(\Omega)}$.

Consider now a bounded linear functional ${\cal L}$ on $H^{1}_{r,a}(\Omega)$
and assume that $\left\Vert {\cal L}\right\Vert=1$.
For all cube $Q\subset \Omega$ of type $(a)$, define
$L^{2}_{0}(Q)=\left\{f\in L^{2}(Q);\ \int_{Q}f(x)dx=0\right\}$. Then,
for all $f\in L^{2}_{0}(Q)$, $\frac f{\left\Vert f\right\Vert_2 \left\vert
Q\right\vert^{1/2}}$ is a type $(a)$ atom, so
that ${\cal L}$ defines a bounded linear functional on $L^{2}_{0}(Q)$. As a
consequence, there exists $b_{Q}\in L^{2}_{0}(Q)$ such that, for all
$f\in L^{2}_{0}(Q)$,
\[
{\cal L}f=\int_{Q}f(x)b_{Q}(x)dx.
\]
Moreover, $\left\Vert b_{Q}\right\Vert_{2}\leq \left\vert
Q\right\vert^{1/2}$. Similarly, if $Q$ is a type $(b)$ cube, ${\cal L}$
defines a bounded linear functional on $L^{2}(Q)$ and there exists
$b_{Q}\in L^{2}(Q)$ with $\left\Vert b_{Q}\right\Vert_{2}\leq \left\vert
Q\right\vert^{1/2}$ such that, for all $f\in L^{2}(Q)$,
\[
{\cal L}f=\int_{Q}f(x)b_{Q}(x)dx.
\]
Observe that, whenever $Q_{1}$ and $Q_{2}$ are type $(b)$ cubes,
$b_{Q_{1}}$ and $b_{Q_{2}}$ coincide on $Q_{1}\cap Q_{2}$. Indeed,
whenever $f\in L^{2}(Q_{1})\cap L^{2}(Q_{2})$, one has
\[
{\cal L}f=\int f(x)b_{Q_{1}}(x)dx=\int f(x)b_{Q_{2}}(x)dx.
\]
Similarly, whenever $Q_{1}$ is a type $(a)$ cube and $Q_{2}$ is a
type $(b)$ cube, $b_{Q_{1}}-b_{Q_{2}}$ is constant on $Q_{1}\cap Q_{2}$.
Indeed, whenever $f$ is supported in $Q_{1}\cap Q_{2}$, $f\in
L^{2}(Q_{1})\cap L^{2}(Q_{2})$ and has zero integral,
one has
\[
{\cal L}f=\int f(x)b_{Q_{1}}(x)dx=\int f(x)b_{Q_{2}}(x)dx.
\]
The key observation at this point is that  for any   $x\in \Omega$,
there exists a type $(b)$ cube that contains $x$ and one defines
$b(x)=b_{Q}(x)$ where
$Q$ is any such cube. This definition is
consistent because of the previous remark. 

Consider now a type $(a)$
cube $Q$. Then, $Q$ is contained in a type $(b)$ cube, hence there
exists
$c_{Q}\in
\complex$ such that $b=b_{Q}+c_{Q}$ on $Q$. One has
\[
\frac 1{\left\vert Q\right\vert}\int_{Q} \left\vert
b(x)-c_{Q}\right\vert^{2}dx =\frac 1{\left\vert Q\right\vert} \int_{Q}
\left\vert b_{Q}(x)\right\vert^2dx \leq 1.
\]
If $Q$ is a type $(b)$ cube,
\[
\frac 1{\left\vert Q\right\vert}\int_{Q} \left\vert
b(x)\right\vert^{2}dx=\frac 1{\left\vert Q\right\vert}\int_{Q} \left\vert
b_{Q}(x)\right\vert^{2}dx  \leq 1.
\]
Hence, $\left\Vert b\right\Vert_{BMO_{z,a}(\Omega)}\leq 1$. One
easily checks that 
\[
{\cal L}f=\int f(x)b(x)dx
\]
whenever $f$ is a finite linear combination of atoms of type $(a)$
or $(b)$ in $H^{1}_{r,a}(\Omega)$.

\paragraph{Assertion $(b)$:} Let $\phi$ be a function in
$BMO_z(\Omega)$. Denote by $D(\reel^n)$ the vector space generated by
$H^1(\reel^n)$-atoms and  by $D_r(\Omega)$ the space of restrictions to
$\Omega$ of functions in $D(\reel^n)$. By definition and density of
$D(\reel^n)$ in
$H^1(\reel^n)$, one has that $D_r(\Omega)$ is dense
in $H^{1}_{r}(\Omega)$ and for
$f\in D_{r}(\Omega)$, $\|f\|_{H^{1}_{r}(\Omega)} = \inf
{\|F\|_{H^1(\reel^n)}}$ where the infimum is taken over all $F \in
D(\reel^n)$ which coincide with $f$ on $\Omega$. For $f\in
D_{r}(\Omega)$, define
\[
{\cal L}f=\int_{\Omega} f(x)\phi(x)dx.
\]
Then, for any function $F\in D(\reel^n)$ which coincides with $f$
on $\Omega$, one has
\[
\left\vert {\cal L}f\right\vert =\left\vert \int_{\reel^n}
F(x)\phi(x)dx\right\vert \leq \left\Vert
F\right\Vert_{H^{1}(\reel^n)} \left\Vert
\phi\right\Vert_{BMO_{z}(\Omega)},
\]
which shows that
\[
\left\vert {\cal L}f\right\vert  \leq
\|f\|_{H^{1}_{r}(\Omega)} \left\Vert
\phi\right\Vert_{BMO_{z}(\Omega)}.
\]
Hence the dual of $H^1_{r}(\Omega)$ contains $BMO_{z}(\Omega)$.

Conversely, let ${\cal L}$ be a bounded linear functional in
$H^{1}_{r}(\Omega)$. For all $f\in H^{1}(\reel^n)$, define
\[
\widetilde{{\cal L}}(f)={\cal L}(\left. f\right\vert_{\Omega}).
\]
The definition of the norm in $H^{1}_{r}(\Omega)$ shows that
$\widetilde{{\cal L}}$ is a bounded linear functional on 
$H^{1}(\reel^n)$:  for all $f\in H^{1}(\reel^n)$,
\[
\left\vert \widetilde{{\cal L}}(f)\right\vert = \left\vert {\cal L}(\left.
f\right\vert_{\Omega}) \right\vert \leq \left\Vert {\cal L}\right\Vert
\left\Vert \left. f\right\vert_{\Omega}\right\Vert_{H^{1}_{r}(\Omega)}
\leq \left\Vert {\cal L}\right\Vert
\left\Vert  f\right\Vert_{H^{1}(\reel^n)}.
\]
Therefore, $\left\Vert \widetilde{{\cal L}} \right\Vert \leq  \left\Vert
{\cal L}\right\Vert$.
Since $BMO(\reel^n)$ is the dual of $H^1(\reel^n)$, there exists
$\phi\in BMO(\reel^n)$ such that, for any finite linear combination of
$H^{1}(\reel^n)$-atoms,
\[
\widetilde{{\cal L}}(f)=\int f(x)\phi(x)dx.
\]
 As a consequence,
\[
\left\Vert \phi\right\Vert_{BMO(\reel^n)}\leq \left\Vert {\cal L}
\right\Vert.
\]
Observe that, if $f$ is any $H^{1}(\reel^n)$-atom supported outside
$\Omega$, $\widetilde{{\cal L}}(f)=0$, which shows that $\phi$ is constant
on each connected component of $^c\Omega$. Fix two such components
$C\ne C'$.
 We let $c,c'$ be the
value of $\phi$ on  $C, C'$ respectively.  Let
$Q$ and
$Q'$ be two cubes of same size respectively contained in $C$ and
$C'$.  Define $a(x)=1$, $x \in Q$ and $a(x)=-1$, $x
\in Q'$. Then, $a$ is a multiple of an  $H^1(\reel^n)$-atom with
support contained in the  smallest cube of $\reel^n$ containing $Q$ and
$Q'$. Since its support is contained outside of $\Omega$, we
have
$\widetilde{{\cal L}}(a)=0$. By construction of $a$, we have $ 
\widetilde{{\cal L}}(a)= c|Q| - c'|Q'|$. Hence $c=c'$ and $\phi=c$ 
outside $\Omega$.

Let $\tilde \phi=\phi - c$. Clearly, 
$\tilde\phi_{|_\Omega} \in BMO_z(\Omega)$.  If $f= F_{|_\Omega}$ with
$F\in D(\reel^n)$, one has
\[
{\cal L}(f)=\widetilde{{\cal L}}(F)=\int_{\reel^n}
F(x)\tilde\phi(x)dx=\int_{\Omega} f(x)\tilde\phi(x)dx.
\]
This proves   $\left(H^1_r(\Omega)\right)^{\prime}
\subset BMO_z(\Omega)$.

We now show that the dual of $VMO_z(\Omega)$ is $H^{1}_{r}(\Omega)$.
This is a consequence of the following Banach space principle. If
$X$ is a Banach space and $Y$ is a closed subspace of $X$, then
$Y'$ is isometric to $X'/ Y^\perp$, where $Y^\perp=\{ L \in X': L(y)=0
\ \forall \, y \in Y\}$. Here, we have $X= VMO(\reel^n)$ and 
$Y=VMO_z(\Omega)$. 

\paragraph{Assertion $(c)$:} is in \cite{coifw} with minor changes due
to our modification of definition.

\paragraph{Assertion $(d)$:} We apply the above abstract principle
with 
$X= H^1(\reel^n)$ and $Y=H^1_z(\Omega)$.

\paragraph{Assertion $(e)$:}  That $(H^1_{r,a}(\Omega))' \subset
BMO_{z,a}(\Omega)$ is straightforward. The converse is already observed
in Section \ref{BMO}.

\section{Local Hardy and $BMO$ spaces on strongly Lipschitz domains}
\label{local} 

We now give localized versions of the previous results.

\subsection{Local Hardy spaces}

We first recall the
definition of
$h^{1}(\reel^n)$ and its atomic decomposition from
\cite{goldberg}.

\paragraph{Definition of $h^{1}(\reel^n)$:}   Let
$\phi\in {\cal S}(\reel^n)$ be a function such that
$\int_{\reel^n} \phi(x)dx=1$. For all $t>0$, define $\phi_t(x)=t^{-n}
\phi(x/t)$. A locally integrable function $f$ on $\reel^n$ is said to be
in $h^1(\reel^n)$ if the maximal function
\[
mf(x)=\sup\limits_{0<t<1} \left\vert \phi_t \ast f(x)\right\vert
\]
belongs to $L^1(\reel^n)$. If it is the case, define
\[
\left\Vert f\right\Vert_{h^1(\reel^n)} = \left\Vert mf\right\Vert_1.
\]
One has $H^{1}(\reel^n)\subset h^{1}(\reel^n)$. It should be noted that a
function in $h^{1}(\reel^n)$ does not
necessarily have zero integral. We note that other maximal functions 
$\sup\limits_{0<t<\delta} \left\vert \phi_t \ast f(x)\right\vert$ with
$\delta>0$ would lead to an equivalent norm.

Replacing $t>0$ by $0<t<1$ in (\ref{semigroup}), one obtains a
characterization
of $h^{1}(\reel^n)$ in terms of a non tangential maximal function
associated with the heat or the Poisson semigroup generated by
$\Delta$ (see \cite{goldberg}).

\paragraph{Atomic decomposition  of $h^{1}(\reel^n)$:}  A
function $a$ is an $h^{1}(\reel^n)$-atom if it is supported in a cube
$Q$, satisfies $\left\Vert a\right\Vert_{2}\leq {\left\vert
Q\right\vert^{-1/2}}$ and has mean-value zero if $\ell(Q)<1$. Then,
$f \in h^{1}(\reel^n)$ if and only if
$f=\sum_{Q} \lambda_{Q}a_{Q}$, where the $a_Q$'s are
$h^1(\reel^n)$-atoms and $\sum_Q|\lambda_Q|<\infty$. Moreover, 
$\left\Vert f\right\Vert_{h^{1}(\reel^n)}$ is comparable with the
infimum of
$\sum_Q \left\vert \lambda_{Q}\right\vert$ taken over all
such decompositions.

\bigskip

We now turn to local Hardy spaces on $\Omega$.  As for
global spaces, three categories of local Hardy spaces on $\Omega$ may
be considered. The first category are restriction spaces.

\paragraph{Definition of $h^{1}_{r}(\Omega)$ and $h^{1}_{z}(\Omega)$:} 
The spaces
$h^{1}_{r}(\Omega)$ and
$h^{1}_{z}(\Omega)$ are defined in the same way as $H^{1}_{r}(\Omega)$
and $H^{1}_{z}(\Omega)$, replacing $H^{1}(\reel^n)$ by $h^{1}(\reel^n)$.
Observe that $h^{1}_{z}(\Omega)$ is a strict subspace of
$h^{1}_{r}(\Omega)$ (see \cite{dafni}, Proposition 6.4).
\bigskip

The second category is made up of atomic spaces. Here, we adopt 
 definitions different from \cite{stkrch} and \cite{dafni}. We feel they
are more natural ones.
 
\paragraph{Definition of type $(a)$ and type 
$(b)$ local cubes:}
  Let $\delta>0$.  A cube $Q$ is
a  type $(a)$ local cube if $\ell(Q)< \delta$ and $4Q\subset \Omega$,  a 
type $(b_{far})$ local cube if $\ell(Q)\geq \delta$ and $4Q\subset \Omega$,
and  a type $(b_{close})$ local cube if $2Q\subset \Omega$ and
$4Q\cap 
\partial\Omega \neq \emptyset$. The type $(b)$ local cubes are those of
type $(b_{far})$ or $(b_{close})$. We always arrange
$\delta$ so that the class of type $(b_{far})$ local cubes is not empty.
To simplify the exposition, we fix $\delta=1$.

\paragraph{Definition of type $(a)$ and type 
$(b)$ local atoms:} A measurable function $a$ on $\Omega$ is called a
type
$(a)$ local atom if it is supported in a type $(a)$ local cube $Q$ with
\[
\left\Vert a\right\Vert_2\leq \left\vert Q\right\vert^{-1/2}\mbox{ 
and } \int a(x)dx=0.
\]
A measurable function $a$ on $\Omega$ is called a type $(b_{far})$
(resp. $(b_{close})$)  local atom if it is supported in a type
$(b_{far})$ (resp. $(b_{close})$) local  cube
$Q$ with
\[
\left\Vert a\right\Vert_2\leq \left\vert Q\right\vert^{-1/2}.
\]
Note that   type $(b)$  local atoms  do not have
mean value zero.

\paragraph{Definition of $h^{1}_{r,a}(\Omega)$:} A function $f$
defined on $\Omega$ belongs to
$h^{1}_{r,a}(\Omega)$ if
\[
f=\sum\limits_{(a)} \lambda_{Q}a_{Q} + \sum\limits_{(b)} \mu_{Q}b_{Q}
\]
where the $a_Q$'s are type $(a)$ local atoms, the $b_{Q}$'s are type
$(b)$ local atoms and
$\sum\limits_{(a)}\left\vert \lambda_Q\right\vert +
\sum\limits_{(b)}\left\vert \mu_{Q}\right\vert <+\infty$.
Define $\left\Vert f\right\Vert_{h^1_{r,a}}$ as the infimum of
$\sum\limits_{(a)}\left\vert \lambda_Q\right\vert +
\sum\limits_{(b)}\left\vert \mu_{Q}\right\vert$ over
all such decompositions.

\paragraph{Definition of $h^{1}_{z,a}(\Omega)$:} A function $f$
defined on $\Omega$ belongs to
$h^{1}_{z,a}(\Omega)$ if
\[
f=\sum\limits_{(a)} \lambda_{Q}a_{Q} + \sum\limits_{(b_{far})}
\mu_{Q}b_{Q}
\]
where the $a_Q$'s are type $(a)$ local atoms, the $b_{Q}$'s are type
$(b_{far})$ local atoms and
$\sum\limits_{(a)}\left\vert \lambda_Q\right\vert +
\sum\limits_{(b_{far})}\left\vert \mu_{Q}\right\vert <+\infty$.
Define $\left\Vert f\right\Vert_{h^1_{z,a}}$ as the infimum of
$\sum\limits_{(a)}\left\vert \lambda_Q\right\vert +
\sum\limits_{(b_{far})}\left\vert \mu_{Q}\right\vert$ over
all such decompositions.

\paragraph{Definition of $h^{1}_{CW}(\Omega)$:} An
$h^1_{CW}(\Omega)$-atom is a function $a$  supported in
$Q\cap\overline\Omega$, where $Q$ is a cube centered in $\Omega$
(but not necessarily included in $\Omega$) with
\[
\left\Vert a\right\Vert_{2} \leq \frac
1{{\left\vert Q \cap\Omega\right\vert}^{1/2}},
\qquad
{\rm and}\qquad {\rm if}\ \ell(Q)< 1,
\int a(x)dx=0  .
\]
A function $f$ is in $h^1_{CW}(\Omega)$ if it can be written
\[
f=\sum\limits_{Q} \lambda_Q a_Q.
\]
where the $a_Q$'s are $h^1_{CW}(\Omega)$-atoms and 
$\sum\limits_{Q} |\lambda_Q| <\infty$. The norm is defined as usual.

\begin{remark}\label{remarkbdd}
 Each global  Hardy space is contained in the
corresponding local space. Also, $h^1_{z,a}(\Omega)$ and $h^1_{z}(\Omega)$
are respectively strict subspaces of $h^1_{r,a}(\Omega)$ and
$h^1_{r}(\Omega)$. If
$\Omega$ is bounded, one can see that $h^1_{r,a}(\Omega) =
H^1_{r,a}(\Omega)$, $h^1_{r}(\Omega) = H^1_{r}(\Omega)$ and that 
$h^1_{CW}(\Omega)= H^1_{CW}(\Omega) + \complex {\cal X}_\Omega$,
$h^1_{z,a}(\Omega)= H^1_{z,a}(\Omega) + \complex {\cal X}_\Omega$, and
$h^1_{z}(\Omega)= H^1_{z}(\Omega) + \complex {\cal X}_\Omega$. Here
${\cal X}_\Omega$ is the indicator function of $\Omega$. All these facts
but the  inclusion
$h^1_{z}(\Omega)\subset H^1_{z}(\Omega) + \complex {\cal X}_\Omega$ are
easy to prove.  For the latter one uses the following observation using
maximal characterizations: if
$f\in h^1(\reel^n) $ has compact support and vanishing mean, then $f \in
H^1(\reel^n)$. Details are left to the reader. 
\end{remark}

These local Hardy spaces compare as follows.

\begin{theorem} \label{h1equalities}
\begin{itemize}
\item[$(a)$] $h^1_{r}(\Omega) = h^1_{r,a}(\Omega)$.
\item[$(b1)$]
$h^1_{z,a}(\Omega) = h^1_{CW}(\Omega)$.
\item[$(b2)$]
$h^1_{z,a}(\Omega) =  h^1_{z}(\Omega) $.
\end{itemize}
\end{theorem}

Note that (a) holds with no
restriction on
$\Omega$ while it is not true for global Hardy spaces.

 We admit this
result for the moment but   the inclusion $  h^1_{z}(\Omega) \subset
h^1_{z,a}(\Omega) 
$, which will be seen in the course of proving the next theorem.

\bigskip

Finally, we consider the third category of local Hardy spaces.

\paragraph{Definition of $h^{1}_{\max, L}(\Omega)$:} If
$L=(A,\Omega,V)$ is a second order elliptic operator in divergence
form  and if $f\in L^{1}_{loc}(\Omega)$ with
slow growth, we  say that
$f\in h^{1}_{max,L}(\Omega)$ if
\[
f^{\ast}_{loc,L}(x)=\sup\limits_{\left\vert y-x\right\vert <t\le1}
\left\vert e^{-tL^{1/2}}f(y)\right\vert \in L^{1}(\Omega).
\]
Define
\[
\left\Vert f\right\Vert_{h^{1}_{max,L}}=\left\Vert
f^{\ast}_{loc,L}\right\Vert_{1}.
\]
It is evident that $H^{1}_{max,L}(\Omega) \subset
h^{1}_{max,L}(\Omega)$.

The local version of Theorem \ref{comparison} is as follows.
\begin{theorem} \label{comparloc}
Let $\Omega= \reel^n$ or $\Omega$ be a  strongly Lipschitz
domain, and
$L=(A,\Omega,V)$ satisfying $(G_{\infty})$.
\begin{itemize}
\item[$(a)$] One has
$h^1(\reel^n)=h^1_{max,L}(\reel^n)$.
\item[$(b)$]
Assume that $L$ satisfies the DBC.
Then, one has
$h^1_{r,a}(\Omega)=h^1_{max,L}(\Omega)
=h^1_r(\Omega)$.
\item[$(c)$]
Assume that $L$ satisfies the NBC.
Then, one has
$h^1_{z,a}(\Omega)=h^1_{max,L}(\Omega)=h^{1}_{z}(\Omega)$.
\end{itemize}
\end{theorem}

\begin{proposition}\label{comparlocbdd} For a bounded Lipschitz domain,
statements (b) and (c) in  Theorem \ref{comparloc} hold  for
$L=(A,\Omega,V)$ satisfying
$(G_{1})$. 
\end{proposition}

This result applies when the coefficients of $A$ are complex-valued
BUC functions or in the closure of  BUC in $bmo$ (See \cite{estim}).

\begin{remark} When $\Omega$ is bounded and $L$ satisfies
$(G_{1})$,  $h^{1}_{max,L}(\Omega) =
H^{1}_{max,L}(\Omega)$. Indeed, one inclusion holds. For the converse,  
 consider $x\in \Omega$, $t>1$ and $y\in
\Omega$ satisfying $\left\vert y-x\right\vert<t$. Lemma
\ref{Poissonkernel}  in Appendix A yields
\[
\left\vert P_tf(y)\right\vert \leq \int_{\Omega} \frac{Ct}{t+\left\vert
y-z\right\vert}\left\vert f(z)\right\vert dz \leq C\left\Vert
f\right\Vert_1.
\]
As a consequence, for all $x\in \Omega$,
\[
\left\vert f^{\ast}_L(x)\right\vert \leq C\left(\left\vert
f^{\ast}_{loc,L}(x)\right\vert + \left\Vert f\right\Vert_1\right)
\]
and
\[
\left\Vert f^{\ast}_L\right\Vert_1 \leq C^{\prime}\left(\left\Vert
f^{\ast}_{loc,L}\right\Vert_1+\left\Vert f\right\Vert_1\right).
\]
\end{remark}

The strategy to prove Theorem \ref{comparloc} and Proposition
\ref{comparlocbdd} is essentially the same  as for the global spaces: we
need a few local
$bmo$-spaces and some duality results, comparison between maximal
functions and  area functionals, and the theory of tent spaces.

\begin{remark}   Assertion  $(c)$ in 
Proposition
\ref{comparlocbdd}   applies to the Neumann
Laplacian on a bounded $\Omega$. Together with  Remark \ref{remarkbdd} this
completes the proof of Theorem
\ref{H1equalities}, $(b2)$. 
\end{remark}

\subsection{$bmo$ spaces}\label{bmo}

A locally square-integrable function $f$ on $\reel^n$ is said to be in
$bmo(\reel^n)$ if
\[
\left\Vert
\phi\right\Vert_{bmo(\reel^n)}^2=\sup\left(\sup\limits_{\ell(Q)
\le 1}
\frac 1{\left\vert Q\right\vert} \int_Q \left\vert
\phi(x)-\phi_Q\right\vert^2 dx, 
\sup\limits_{\ell(Q) > 1}
\frac 1{\left\vert Q\right\vert} \int_Q \left\vert
\phi(x)\right\vert^2 dx\right) <+\infty.
\]
Define $vmo(\reel^n)$
as  the closure
of $C_c(\reel^n)$  in
$bmo(\reel^n)$.  
It is well-known that $bmo(\reel^n)$ is the dual of $h^1(\reel^n)$,
which is the dual of $vmo(\reel^n)$ \cite{goldberg}.

\smallskip

 Define $bmo_z(\Omega)$, $vmo_z(\Omega)$ and $bmo_r(\Omega)$
analogously to the corresponding global $BMO$ or $VMO$ spaces,
replacing 
$BMO(\reel^n)$   by $bmo(\reel^n)$ and $VMO(\reel^n)$   by
$vmo(\reel^n)$.

\smallskip

A locally square-integrable function $f$ on $\Omega$ is in 
$bmo_{z,a}(\Omega)$ if
\[
\left\Vert
\phi\right\Vert_{bmo_{z,a}(\Omega)}^2=\sup\left(\sup\limits_{(a)}
\frac 1{\left\vert Q\right\vert} 
\int_{Q} \left\vert \phi(x)-\phi_{Q}\right\vert^2 dx , 
\sup\limits_{(b)}  \frac
1{\left\vert Q\right\vert} 
\int_{Q}\left\vert\phi(x)\right\vert^2 dx \right)<+\infty,
\]
where $\sup\limits_{(a)}$ (resp. $\sup\limits_{(b)}$) means 
that the supremum is taken over all 
 type $(a)$ (resp.  
 $(b)$) local cubes.
 
\smallskip

A  locally square-integrable function $f$ on $\Omega$ is in 
$bmo_{r,a}(\Omega)$ if
\[
\left\Vert \phi\right\Vert_{bmo_{r,a}(\Omega)}^2=\sup\left(
\sup\limits_{(a)}
\frac 1{\left\vert Q\right\vert} 
\int_{Q} \left\vert \phi(x)-\phi_{Q}\right\vert^2 dx 
, 
\sup\limits_{(b_{far})}  \frac
1{\left\vert Q\right\vert} 
\int_{Q}\left\vert\phi\right\vert^2
\right)
<+\infty.
\] 
\smallskip

A locally square-integrable function $\phi$ defined on $\Omega$ is in
$bmo_{CW}(\Omega)$ if
\[
\left \Vert \phi\right\Vert_{bmo_{CW}(\Omega)}^2 =
\sup\left( 
\sup_{\ell(Q)<1} \frac 1{\left\vert Q\cap\Omega\right\vert}
\int_{Q\cap\Omega}
\left\vert \phi(x)-\phi_{Q\cap\Omega}\right\vert^2 dx
, 
\sup\limits_{{\ell(Q)\ge 1}}  \frac 1{\left\vert Q\cap\Omega\right\vert}
\int_{Q\cap\Omega}
\left\vert\phi\right\vert^2
\right)<+\infty,
\]
where the  cubes  have center in
$\Omega$.  The space $vmo_{CW}(\Omega)$ is defined as  the closure of
$C_c(\Omega)$ in $bmo_{CW}(\Omega)$.
\smallskip

The duality results for local spaces and the comparisons between $bmo$
spaces are the same as for the global spaces.  Let us state them for
completeness.

\begin{theorem} \label{thdualloc} 
\begin{itemize}
\item[$(a)$]
The dual of $h^1_{r,a}(\Omega)$ is $bmo_{z,a}(\Omega)$.
\item[$(b)$]
The dual of $h^1_{r}(\Omega)$ is $bmo_{z}(\Omega)$, the dual of 
$vmo_{z}(\Omega)$ is $h^{1}_{r}(\Omega)$.
\item[$(c)$]
The dual of $h^{1}_{CW}(\Omega)$ is $bmo_{CW}(\Omega)$, the dual of
$vmo_{CW}(\Omega)$ is $h^{1}_{CW}(\Omega)$. 
\item[$(d)$]
The dual of $h^1_z(\Omega)$ is $bmo_r(\Omega)$.
\item[$(e)$]
The dual of $h^1_{z,a}(\Omega)$ is $bmo_{r,a}(\Omega)$.
\end{itemize}
\end{theorem}

\begin{theorem} \label{bmoequalities}
\begin{itemize}
\item[$(a)$]
$bmo_{z,a}(\Omega) = bmo_{z}(\Omega)$. 
\item[$(b1)$]
$bmo_{CW}(\Omega)=  bmo_{r,a}(\Omega) $.
\item[$(b2)$]
$bmo_{r}(\Omega) =   bmo_{r,a}(\Omega) $.
\end{itemize}\smallskip
\end{theorem}

Again, the difficult parts are $(h^1_{z,a}(\Omega))' \subset
bmo_{r,a}(\Omega)$ and $   bmo_{r,a}(\Omega) \subset bmo_{r}(\Omega)
$, which are proved using Theorem \ref{thdualloc}, $(c)$ and 
Theorem \ref{bmoequalities}, $(b1)$.

\subsection{Proofs of equalities between local Hardy spaces}

\paragraph{Proof of Theorem \ref{comparloc}.}

In each case, the most involved part is to imbed our maximal space into
an atomic space. We concentrate on this.

One has the local statements corresponding to the results in Sections 
\ref{area} and \ref{BMOCarleson}, in which  the square functions and the
Carleson measures are truncated at some fixed time $t<t_0$, say for
example
$t=1$. Except for some technical adjustments the proofs are the same
and left to the reader. 

The idea is to use this in the representation formula  $f=f_1+f_2$ where
\[
{f_1}=4\int_0^{1} (tL^{1/2} P_t)(tL^{1/2} P_tf)\frac{dt}t 
\]
and $f_2=  \frac 14
(2L^{1/2}+I)P_1 P_1f$
  for $f \in h^1_{max,L}(\Omega) \cap
L^2(\Omega)$.\footnote{If $\Omega$ is bounded and $V=W^{1,2}(\Omega)$
then the formula holds if $\int_\Omega f =0$. If the mean of $f$ is
not zero, then it applies to $\tilde f = f - c{\cal X}_\Omega$ with
the constant $c$ so that the mean of $\tilde f$ is zero. Conclude
with  ${\cal X}_\Omega \in h^1_{CW}(\Omega)$.}
 
 For $f_1$, we proceed using the tent spaces again and then
eliminate the requirement that $f \in L^2(\Omega)$ to obtain 
$f_1 \in h^1(\reel^n)$ or $h^1_{r,a}(\Omega)$ under DBC or
$h^1_{CW}(\Omega)$ under NBC. 

 Let us consider
$f_2$. The idea is to prove  that $f_2 \in h^1_{CW}(\Omega)$ in
each case. Indeed, when
$\Omega=\reel^n$ we have $h^1_{CW}(\reel^n) = h^1(\reel^n)$, under
DBC $h^1_{CW}(\Omega) \subset  h^1_{r,a}(\Omega)$  from the
definitions of theses spaces and easy arguments,  and under
NBC, 
$h^1_{CW}(\Omega) =  h^1_{z,a}(\Omega)$  from Theorem
\ref{h1equalities}, (b1).

Here is the argument. Assume first  
$\Omega=\reel^n$. Observe that $P_1f=g$ is bounded by $f^*_{loc,L}$
which is in $L^1(\reel^n)$. Also since $L^{1/2}P_1= -
\partial_tP_t\mid_{t=1}$ the subordination formula yields that the
kernel $K(x,y)$  of $(2L^{1/2}+I)P_1$ is bounded by 
$ck(x,y)$ with $k(x,y)=(1+|x-y|)^{-n-1}$.

Take $(Q_k)$ be a covering of $\reel^n$ by cubes with size 1 obtained
by translation from the unit cube $[0,1]^n$. Let $(\eta_k)$ be a smooth
partition of unity associated with this covering so that $\eta_k$ is
supported in $2Q_k$. Then one has
$$
f_2(x) = \sum_{k} b_{k}(x)
$$
where
$$
b_{k}(x) = \eta_k(x) \int_{\reel^n} K(x,y)  g(y)\, dy
$$
Observe that $b_{k}$ is supported in $2Q_k$.
Set 
$\lambda_{k} = |2Q_k|^{1/2} \left(\int |b_{k}|^{2}\right)^{1/2}.
$
We have 
$$
\lambda_{k} \le c |2Q_k|^{1/2} \left(\int_{2Q_k}  |\int_{\reel^n}
k(x,y)
| g(y)|\, dy|^2 dx\right)^{1/2}.
$$
Observe that
$k(x,y) \le c \inf_{|x-z| \le 2} k(z,y)
$
for all $x, y \in \reel^n$. For $x\in 2Q_k$, since $\ell(Q_k) \le 2$, 
$k(x,y) \le c \inf_{z \in 2Q_k} k(z,y)
$
for all $y$.  Hence 
$$
\lambda_{k} \le c |2Q_k| \int_{\reel^n}
\inf_{z \in 2Q_k} k(z,y)
| g(y)|\, dy \le \iint {\cal X}_{2Q_k}(x) 
k(x,y) |g(y)|
\, dydx
$$
and it follows from the finite overlap property of the family $(2Q_k)$
that 
$$
\sum_{k}  \lambda_{k} \le  c \iint k(x,y) |g(y)| \, dydx \le 
c\|g\|_1 \le c\|f\|_{h^1_{max,L}(\reel^n)}.$$

If we set $a_k =\lambda^{-1}_k b_k$ when $\lambda_k\ne 0$, then $a_k$
is an $h^1(\reel^n)$-atom. Thus, $f_2 \in h^1(\reel^n)$ with
$\|f_2\|_{h^1(\reel^n)} \le c\|f\|_{h^1_{max,L}(\reel^n)}$.

Assume now that $\Omega$ is strongly Lipschitz and $L$ satisfies 
 either boundary condition. Let  $(Q_k)$ be  a covering of
$\Omega$ with cubes of $\reel^n$ such that $\ell(Q_k) = 1$.
We keep only those cubes which intersect $\Omega$. If $Q_k$ has
center in $\Omega$, we set $\lambda_k=|2Q_k\cap \Omega|^{1/2} \left(\int
|b_{k}|^{2}\right)^{1/2}$. If $Q_k$ has center outside of $\Omega$,
then   we replace $Q_k$ by $\tilde Q_k$ with  center in $Q_k \cap
\Omega$ and  $\ell(\tilde Q_k)= 2\ell(Q_k)$ and define
$\lambda_k=|2\tilde Q_k\cap \Omega|^{1/2} \left(\int
|b_{k}|^{2}\right)^{1/2}$. With these modifications of $\lambda_k$, we
see from the same argument that  $a_k$ is an $h^1_{CW}(\Omega)$-atom and
that
$\sum |\lambda_k| \le c\|f\|_{h^1_{max,L}(\Omega)}$ remarking that all
integrals should take place on $\Omega$.

\paragraph{Proof of Proposition \ref{comparlocbdd}.} 

We want to relax the condition $(G_\infty)$ to
$(G_1)$ when $\Omega$ is bounded. The same arguments works once we make
sure of small time decay estimates for the Poisson kernel. This is
proved in  Appendix A.

\paragraph{Proof of Theorem \ref{h1equalities}.} 

The proof of $(a)$ is as in the global case and is skipped.

The proofs of $h^1_{z,a}(\Omega) \subset h^1_{CW}(\Omega)$ and 
$h^1_{z,a}(\Omega)\subset h^1_{z}(\Omega)$ are straightforward from the
definitions.

It remains to prove $h^1_{CW}(\Omega) \subset  h^1_{z,a}(\Omega)$.
Because of Remark \ref{remarkbdd} and the global case, this is already
known if $\Omega$ is bounded. We assume next that $\Omega$ is unbounded.
 Let
$a$ be an
$h^1_{CW}$-atom supported on a cube
$Q$ centered in $\Omega$. Since $\Omega$ is unbounded, we have $|Q\cap \Omega|
\sim |Q|$ (see Appendix B).

If $\ell(Q) < 1$, we proceed as in the global case. Either $a$ is a
type $(a)$ local atom or can be decomposed into a sum of type $(a)$
local atoms.

Assume $\ell(Q)\ge 1$. If $a$ is a type $(b_{far})$ local atom, we are
done. It remains to argue when $Q$ is close to the boundary, ie $4Q \cap
\partial \Omega \ne \emptyset$. In this case, we claim there exists a
type $(b_{far})$ local cube $Q'$ with $\ell(Q') =\ell(Q)$ and the
distance between $Q$ and $Q'$ is comparable to $\ell(Q)$ (See Appendix
B). Define $\tilde a = a - (\frac 1{|Q'|} \int_{Q\cap\Omega} a) {\cal
X}_{Q'}$. Then $\tilde a$ is a multiple of $H^1_{CW}(\Omega)$-atom
(supported in
$cQ \cap\Omega$ for some constant $c$ that does not depend on $Q$),
thus $\tilde a \in H^1_{z,a}{\Omega} \subset h^1_{z,a}{\Omega}$.
Now, $\tilde a  -a$ clearly is a type $(b_{far})$ local  atom, hence it
belongs to 
$h^1_{z,a}({\Omega})$.

\section{Other maximal functions} \label{other}
 As a consequence of the atomic decomposition for the maximal space, one
may use other maximal
functions, such as the vertical
and the non-tangential maximal functions associated with $e^{-tL}$. More
precisely, the 
following holds:
\begin{theorem} \label{heathardy}
Let $L=(A,\Omega,V)$. Assume that $\Omega=\reel^n$ or that $\Omega$ be a
strongly Lipschitz domain of 
$\reel^n$  under DBC with $^c\Omega$ unbounded or under NBC. Assume also
that
$L$ and $L^*$ satisfy 
$(G_{\infty})$. 
 The following are equivalent:
\begin{equation} \label{first}
\sup\limits_{t>0} \left\vert
e^{-tL}f(x)\right\vert \in L^{1}(\Omega),
\end{equation}
\begin{equation} \label{second}
\sup\limits_{\left\vert x-y\right\vert <\sqrt{t}} \left\vert
e^{-tL}f(y)\right\vert \in L^{1}(\Omega),
\end{equation}
\begin{equation} \label{third}
 f \in H^{1}_{\max,L}(\Omega),
\end{equation}
One also has the analogous local statement, replacing
 $H^{1}_{\max,L}$ by $h^{1}_{\max,L}$ and $t>0$ by
$0<t<t_0$ for any $t_0>0$ without restriction on $\Omega$. Moreover, if
$\Omega$ is bounded then 
$(G_{1})$ suffices.
\end{theorem}

The new assumption that $L^*$ satisfies $(G_\tau)$ simply
means that   (\ref{Holder}) holds for $K_t(y,x)$ too. Again, when
$L$ is real, this is not a supplementary hypothesis.  
 
We write the proof  for global spaces, under DBC 
when $^c\Omega$
is unbounded  for example. 
We have already proved (Theorem \ref{comparison}) the
implication $(\ref{third}) \Rightarrow f\in H^{1}_{r,a}(\Omega)$ and the
implication $f\in H^{1}_{r,a}(\Omega) \Rightarrow (\ref{first})$ is an
easy consequence of the estimates for $K_{t}$ (as the proof of $f\in
H^{1}_{r,a}(\Omega) \Rightarrow f \in H^{1}_{\max,L}(\Omega)$). We
therefore turn to the proofs of $(\ref{first}) \Rightarrow
(\ref{second})$ and  $(\ref{second}) \Rightarrow
(\ref{third})$. 

The argument for $(\ref{first}) \Rightarrow
(\ref{second})$ relies upon
the comparison between the $L^{1}$ norms of two
maximal functions and is inspired by \cite{feffstein}, p.185. 
For all $\alpha>0$ and $v:\Omega\times]0,+\infty[ \to \CC$ set
\[
v^{\ast}_{\alpha}(x)=\sup\limits_{\left\vert y-x\right\vert
<\alpha \sqrt{t}} \left\vert v(y,t)\right\vert.
\]
If $f\in
L^{1}_{loc}$ with slow growth, set
\[
u(x,t)=e^{-tL}f(x),\ u^{+}(x)=\sup\limits_{t>0}\left\vert
u(x,t)\right\vert,\ u^{\ast}(x)=u^{\ast}_{1}(x).
\]
Recall that we assume $(G_\infty)$ so that slow growth insures that $u$
is well-defined.

Finally, for all $\varepsilon>0$, all $N\in \nat$ and all $x\in 
\Omega$, consider
\[
u^{\ast}_{\varepsilon,N}(x)=\sup\limits_{\left\vert y-x\right\vert
<\sqrt{t}<\varepsilon^{-1}} \left\vert u(y,t)\right\vert
\left(\frac{\sqrt{t}}{\sqrt{t}+\varepsilon}\right)^{N}\left(1+\varepsilon
\left\vert y\right\vert\right)^{-N}
\]
and
\[
U^{\ast}_{\varepsilon,N}(x)=\sup\limits_{\left\vert y-x\right\vert
<\sqrt{t}<\varepsilon^{-1},\ \left\vert y^{\prime}-x\right\vert
<\sqrt{t}<\varepsilon^{-1}} \left(\frac{\sqrt{t}}{\left\vert
y-y^{\prime}\right\vert}\right)^\mu \left\vert
u(y,t)-u(y^{\prime},t)\right\vert
\left(\frac{\sqrt{t}}{\sqrt{t}+\varepsilon}\right)^{N}\left(1+\varepsilon
\left\vert y\right\vert\right)^{-N}.
\]
for some $\mu>0$ to be chosen later.

We intend to show the following proposition:
\begin{proposition} \label{pointwise}
There exists $C>0$ such that, for all $f\in L^{1}_{loc}$, $\left\Vert
u^{\ast}\right\Vert_{1}\leq C\left\Vert u^{+}\right\Vert_{1}$.
\end{proposition}
Notice first that the $L^{1}$-norm of $u^{\ast}_{\alpha}$ is
controlled by the $L^{1}$-norm of $u^{\ast}$. More
precisely, the following holds (see \cite{feffstein}, Lemma 1, p. 166):
\begin{lemma} \label{apert}
There exists $C$ such that, for all continuous function $v$ on $\Omega 
\times \left]0,+\infty\right[$ and all $\alpha>0$,
\[
\left\Vert v^{\ast}_{\alpha} \right\Vert_{1}\leq C\alpha^n \left\Vert
v^{\ast}\right\Vert_{1}.
\]
\end{lemma}

Note that this inequality holds if $v$ is truncated for $t>t_0$. 

The proof of Proposition \ref{pointwise} relies on the following observation:
\begin{lemma} \label{gradmax}
Assume that $u_{\varepsilon,N}^{\ast}\in L^{1}$. Then
\[
\left\Vert U_{\varepsilon,N}^{\ast}\right\Vert_{1} \leq C\left\Vert
u_{\varepsilon,N}^{\ast}\right\Vert_{1},
\]
where $C$ is independent on $\varepsilon,N$ and $u$.
\end{lemma}
Fix $x\in \Omega$ and consider $y,y^{\prime}$ and $t$ such that $\left\vert
y-x\right\vert <\sqrt{t}$ and $\left\vert
y^{\prime}-x\right\vert <\sqrt{t}$. Define also
\[
v(y,t)=u(y,t) \left(1+\varepsilon \left\vert 
y\right\vert\right)^{-N}\left(\frac{\sqrt{t}}{\sqrt{t}+\varepsilon}\right)^{N}
{\cal X}_{]0,1[}(\varepsilon\sqrt t)
\]
so that $v^*_1= u_{\varepsilon,N}$.
Start from
\[
u(y,t)-u(y^{\prime},t)=\int_\Omega
\left(K_{t/2}(y,z)-K_{t/2}(y^{\prime},z)\right)u(z,t/2)dz= I_0+
\sum_{k\ge 1} I_k,
\]
where
$$
I_0 =  \int_{\left\vert z-y\right\vert \leq \sqrt{t}}  \left\vert
K_{t/2}(y,z)-K_{t/2}(y^{\prime},z)\right\vert \left\vert
u(z,t/2)\right\vert dz,
$$
and 
$$
I_k = \int_{2^{k-1}\sqrt{t}<\left\vert
z-y\right\vert \leq 2^{k}\sqrt{t}}  \left\vert
K_{t/2}(y,z)-K_{t/2}(y^{\prime},z)\right\vert \left\vert
u(z,t/2)\right\vert dz.
$$
Using $(1+ \varepsilon|z|)^N \le (1+ \varepsilon|y|)^N (1+2^k)^N$
if $\left\vert
z-y\right\vert \leq 2^{k}\sqrt{t}$ and $\varepsilon\sqrt t<1$
and 
$$
\int_{2^{k-1}\sqrt{t}<\left\vert
z-y\right\vert \leq 2^{k}\sqrt{t}}  \left\vert
K_{t/2}(y,z)-K_{t/2}(y^{\prime},z)\right\vert  dz
\le c \left(\frac{\left\vert
y-y^{\prime}\right\vert}{\sqrt{t}}\right)^{\mu} e^{-\alpha
2^{2k}}
$$
for some $\mu>0$ and $\alpha>0$ from $(G_\infty)$ for $L^*$, we easily
get
$$
\left(\frac{\sqrt{t}}{\sqrt{t}+\varepsilon}\right)^{N}
\frac {\left\vert u(y,t)-u(y^{\prime},t)
\right\vert}{\left(1+\varepsilon
\left\vert  y\right\vert\right)^{N}}
\le c \left(\frac{\left\vert
y-y^{\prime}\right\vert}{\sqrt{t}}\right)^{\mu}
 \left(v^{\ast}_{2}(x)+\sum\limits_{k\geq 1} e^{-\alpha
2^{2k}} \left(1+2^{k}\right)^{N}v^{\ast}_{2^{k}+1}(x)\right).
$$
Therefore,
\[
U_{\varepsilon,N}^{\ast}(x)\leq c
 \left(v_{2}^{\ast}(x)+\sum\limits_{k\geq 1} e^{-\alpha
2^{2k}} \left(1+2^{k}\right)^{N}v^{\ast}_{2^{k}+1}(x)\right).
\]
Lemma \ref{gradmax} follows at once from Lemma \ref{apert}. \rule{2mm}{2mm}

We now prove Proposition \ref{pointwise}, following \cite{feffstein}, p. 186. Consider $f$ such that $u^{+}\in
L^{1}$ and $N\in \nat$ large enough, so
that one easily derives that $u^{\ast}_{\varepsilon,N}\in L^{1}$ for all
$\varepsilon>0$. Define $G_{\varepsilon,N}=\left\{x\in
\Omega;U^{\ast}_{\varepsilon,N}(x)\leq Bu^{\ast}_{\varepsilon,N}(x)\right\}$
for some $B>0$ to be chosen. Then, one has
\[
\begin{array}{lll}
\displaystyle \int_{\Omega\setminus G_{\varepsilon,N}}
u^{\ast}_{\varepsilon,N}(x)dx & \leq & \displaystyle \frac 1B
\int_{\Omega\setminus G_{\varepsilon,N}}
U^{\ast}_{\varepsilon,N}(x)dx \\
\\
& \leq & \displaystyle \frac 12 \int_{\Omega} u^{\ast}_{\varepsilon,N}(x)dx
\end{array}
\]
provided that $B$ is large enough.

Moreover, for almost all $x\in G_{\varepsilon,N}$, one has
$u^{\ast}_{\varepsilon,N}(x)\leq CM(x)$, where
\[ 
M(x)=\sup\limits_{Q\ni x} \left(\frac 1{\left\vert Q\cap \Omega\right\vert}
\int_{Q\cap \Omega} u^{+}(y)^r dy\right)^{1/r},
\]
with $0<r<1$ (in this definition, the cubes are centered in
$\Omega$). Indeed, let $x\in G_{\varepsilon,N}$ for which
$u^{\ast}_{\varepsilon,N}(x)<\infty$. There exist
$y,t$ such that $\left\vert y-x\right\vert <\sqrt{t}<\varepsilon^{-1}$
and
\[
\left\vert u(y,t)\right\vert
\left(\frac{\sqrt{t}}{\sqrt{t}+\varepsilon}\right)^{N}
\left(1+\varepsilon\left\vert y\right\vert\right)^{-N} \geq \frac 12
u^{\ast}_{\varepsilon,N}(x).
\]
Since $x\in G_{\varepsilon,N}$, if $\left\vert z-x\right\vert <\sqrt{t}$ and $\left\vert
z^{\prime}-x\right\vert <\sqrt{t}$, one has
\[
\left(\frac {\sqrt{t}}{\left\vert z-z^{\prime}\right\vert}\right)^\mu
\left\vert u(z,t)-u(z^{\prime},t)\right\vert
\left(\frac{\sqrt{t}}{\sqrt{t}+\varepsilon}\right)^{N}
\left(1+\varepsilon\left\vert z\right\vert\right)^{-N} \leq 2B
\left\vert u(y,t)\right\vert
\left(\frac{\sqrt{t}}{\sqrt{t}+\varepsilon}\right)^{N}
\left(1+\varepsilon\left\vert y\right\vert\right)^{-N}
\]
hence,
\[
\left(\frac {\sqrt{t}}{\left\vert z-z^{\prime}\right\vert}\right)^\mu
\left\vert u(z,t)-u(z^{\prime},t)\right\vert \leq c\left\vert
u(y,t)\right\vert.
\]
 It follows that
\[
\left\vert u(z,t)\right\vert \geq \frac 12 \left\vert
u(y,t)\right\vert
\]
when $z\in A=\left\{w;\ \left\vert w-x\right\vert <\sqrt{t}\mbox{
and }\left\vert
w-y\right\vert <\frac {\sqrt{t}}{2C}\right\}$. Therefore, when $z\in
A$, one has
\[
\left\vert u(z,t)\right\vert \geq \frac 12 \left\vert u(y,t)
\right\vert \left(\frac{\sqrt{t}}{\sqrt{t}+\varepsilon}\right)^{N}
\left(1+\varepsilon\left\vert y\right\vert\right)^{-N} \geq \frac 14
u^{\ast}_{\varepsilon,N}(x).
\]
Hence,
\[
\begin{array}{lll}
\displaystyle M(x)^r & \geq & \displaystyle  \frac c{\left\vert
B(x,2\sqrt{t})\right\vert} \int_{B(x,2\sqrt{t})} u^{+}(z)^r dz \\
\\
& \geq & \displaystyle \frac c{\left\vert
B(x,2\sqrt{t})\right\vert} \int_{B(x,2\sqrt{t})} |u(z,t)|^r dz \\
\\
& \geq & \displaystyle c\left(\frac 14
u^{\ast}_{\varepsilon,N}(x)\right)^r \frac{\left\vert A\right\vert}{\left\vert
B(x,2\sqrt{t})\right\vert} \\
\\
& \geq &\displaystyle c u^{\ast}_{\varepsilon,N}(x)^r.
\end{array}
\]
Finally, using the fact that $1/r>1$, one obtains that
\[
\begin{array}{lll}
\displaystyle \int_{\Omega} u^{\ast}_{\varepsilon,N}(x) dx &
\leq & \displaystyle 2 \int_{G_{\varepsilon,N}}
u^{\ast}_{\varepsilon,N}(x) dx \\
\\
& \leq & \displaystyle C \int_{G_{\varepsilon,N}}
M(x) dx \\
\\
& \leq & \displaystyle C \int_{\Omega}
M(x) dx \\
\\
& \leq & \displaystyle C \int_{\Omega}
u^{+}(x) dx
\end{array}
\]
where $C$ does not depend on $\varepsilon$. Letting
$\varepsilon\rightarrow 0$ yields $\|u^\ast\|_1 \le C \|u^+\|_1
$ and $(\ref{second}) $ is proved. 

To complete the proof of Proposition \ref{pointwise}, it
remains to see $(\ref{second})
\Rightarrow (\ref{third})$. This  follows easily from the subordination
formula (\ref{subor}) and Lemma \ref{apert}.  
\rule{2mm}{2mm}

\section*{Appendix A: Kernel estimates} \label{est}
\addcontentsline{toc}{section}{Appendix A:  Kernel estimates}

In this Appendix, we derive some consequences of the Gaussian upper
bounds (\ref{Gauss})  which we assume to hold for $0<t<\tau$. The first 
consequence is that an estimate of the form (\ref{Gauss}) holds for
$t\partial_tK_t(x,y)$ by analyticity of the semigroup (See 
\cite{asterisque}, Chapter I, Lemma 19).

 We first claim the
following:
\begin{lemma} \label{Gaussbigtime}
Assume that $\tau=1$.
Then, for all $t>1$ and all $x,y\in \Omega$, one has
\[
\left\vert K_{t}(x,y)\right\vert \leq Ce^{-\alpha \frac{\left\vert
x-y\right\vert^{2}}t}.
\]
\end{lemma}
The proof relies on the following $L^{2}$-maximum principle (see
\cite{grigo}):
\begin{proposition} \label{maxprinc} Assume that $A \in {\cal A}(c)$.
Let $u(x,t)$ be a function on $\Omega\times \left]0,+\infty\right[$ satisfying
$\partial_{t}u(x,t)+Lu(x,t)=0$ on $\Omega$. Then, if $\xi:\Omega \times
\left]0,+\infty\right[\rightarrow \reel$ is locally Lipschitz and
satisfies the relation
\[
\partial_{t}\xi(x,t)+\alpha \left\vert \nabla \xi(x,t)\right\vert^{2}\leq 0,
\]
where $\alpha=\frac 1{2c^{3}}$, the function
\[
I(t)=\int_{\Omega} \left\vert u(x,t)\right\vert^{2} e^{\xi(x,t)} dx
\]
is non increasing in $t>0$. 
\end{proposition}
Indeed, for all $t>0$, one has
\[
\begin{array}{lll}
I^{\prime}(t) &= & \displaystyle 2\mbox{Re }\int_{\Omega} \partial_{t}u(x,t)
\overline{u(x,t)} e^{\xi(x,t)} dx +\int_{\Omega}
u(x,t) \overline{u(x,t)} \partial_{t}\xi(x,t)e^{\xi(x,t)}dx \\
\\
& = & \displaystyle -2\mbox{Re } \int_{\Omega}Lu(x,t)
\overline{u(x,t)} e^{\xi(x,t)} dx - \alpha \int_{\Omega}
u(x,t) \overline{u(x,t)} \left\vert \nabla_{x}
\xi(x,t)\right\vert^{2} e^{\xi(x,t)}dx \\
\\
& = & \displaystyle -2\mbox{Re } \int_{\Omega} A(x)\nabla u(x,t)
\overline{\nabla u(x,t)} e^{\xi(x,t)} dx - 2 \mbox{Re }
\int_{\Omega} A(x)\nabla u(x,t) \nabla \xi(x,t) \overline{u(x,t)}
e^{\xi(x,t)} dx  \\
\\
& - & \displaystyle \alpha \int_{\Omega}
u(x,t) \overline{u(x,t)} \left\vert \nabla_{x}
\xi(x,t)\right\vert^{2} e^{\xi(x,t)}dx \\
\\
& \leq & \displaystyle -2c \int_{\Omega} \left\vert \nabla
u(x,t)\right\vert^{2} e^{\xi(x,t)} dx +2c^{-1}
\int_{\Omega} \left\vert \nabla u(x,t)\right\vert \left\vert \nabla
\xi(x,t)\right\vert \left\vert u(x,t)\right\vert e^{\xi(x,t)} dx \\
\\
& - & \displaystyle \alpha \int_{\Omega}
\left\vert u(x,t)\right\vert^{2} \left\vert \nabla_{x}
\xi(x,t)\right\vert^{2} e^{\xi(x,t)}dx \\
\\
& \leq & \displaystyle (-2c+c^{-1}/\varepsilon) \int_{\Omega} \left\vert
\nabla
u(x,t)\right\vert^{2} e^{\xi(x,t)} dx +  (\varepsilon c^{-1}-\alpha)
\int_{\Omega}
\left\vert u(x,t)\right\vert^{2} \left\vert \nabla_{x}
\xi(x,t)\right\vert^{2} e^{\xi(x,t)}dx \\
\\
& = & 0.
\end{array}
\]
In the previous computation, $\varepsilon=\alpha c=\frac 1{2c^{2}}$.
Proposition \ref{maxprinc} is therefore proved. \rule{2mm}{2mm}\medskip

In order to prove Lemma \ref{Gaussbigtime}, observe that, for $\beta>0$
small enough (namely, $\beta \leq \frac 1{4\alpha}$), the function
$\xi(x,t)=\beta \frac{\left\vert
x-y\right\vert^{2}}t$ satisfies the assumptions of Proposition
\ref{maxprinc}. Then, write
\[
\begin{array}{lll}
\left\vert K_{t}(x,y)\right\vert & = & \displaystyle \left\vert
\int_{\Omega} K_{t/2}(x,z)e^{\frac{\alpha}2 \frac{\left\vert
x-z\right\vert^{2}}t}
K_{t/2}(z,y) e^{\frac{\alpha}2 \frac{\left\vert z-y\right\vert^{2}}t}
e^{-\frac{\alpha}2 \frac{\left\vert x-z\right\vert^{2}+ \left\vert
z-y\right\vert^{2}}t}dz \right\vert \\
\\
& \leq & \displaystyle \left( \int_{\Omega} \left\vert
K_{t/2}(x,z)\right\vert^{2}e^{\alpha \frac{\left\vert
x-z\right\vert^{2}}t} dz\right)^{1/2} \left( \int_{\Omega} \left\vert
K_{t/2}(z,y)\right\vert^{2}e^{\alpha \frac{\left\vert
z-y\right\vert^{2}}t} dz\right)^{1/2} e^{-\frac{\alpha}4 \frac{\left\vert
x-y\right\vert^{2}}t} \\
& \leq & \displaystyle \left( \int_{\Omega} \left\vert
K_{1/2}(x,z)\right\vert^{2} e^{\alpha \left\vert
x-z\right\vert^{2}} dz\right)^{1/2} \left( \int_{\Omega} \left\vert
K_{1/2}(z,y)\right\vert^{2}e^{\alpha \left\vert
z-y\right\vert^{2}} dz\right)^{1/2} e^{-\frac{\alpha}4 \frac{\left\vert
x-y\right\vert^{2}}t}\\
\\
& \leq & \displaystyle C e^{-\frac{\alpha}4 \frac{\left\vert
x-y\right\vert^{2}}t}.
\end{array}
\]
As a consequence of the upper bounds for $K_t$, we get the following
estimates for the Poisson
kernel:
\begin{lemma} \label{Poissonkernel}
\begin{itemize}
\item[$(a)$]
Assume that $\tau=+\infty$. Then, for all $t>0$, all $x,y\in \Omega$,
\[
\left\vert p_t(x,y)\right\vert \leq \frac {Ct}{\left(t+\left\vert
x-y\right\vert\right)^{n+1}}.
\]
\item[$(b)$]
Assume now that $\tau=1$ and $\Omega$ is bounded. Then, for all $0<t<1$,
all $x,y\in
\Omega$,
\[
\left\vert p_{t}(x,y)\right\vert \leq \frac {Ct}{\left(t+\left\vert
x-y\right\vert\right)^{n+1}}.
\]
For all $t>1$, all $x,y\in \Omega$,
\[
\left\vert p_{t}(x,y)\right\vert \leq \frac {Ct}{\left(t+\left\vert
x-y\right\vert\right)}.
\]
\end{itemize}
By analyticity, the same estimates hold for $t\partial_{t}p_{t}(x,y)$.
\end{lemma}

Just use the subordination formula (\ref{subor}) and the upper estimates for $\left\vert
K_{t}(x,y)\right\vert$.

We now summarize $L^2$-estimates for $\nabla K_t(x,y)$ that follow from
the assumption (\ref{Gauss}) and the Caccioppoli
inequality
(see \cite{autchdom}, Proposition 15):
\begin{proposition} \label{gradient}

\begin{itemize}

\item[$(a)$]
Assume that $\tau=+\infty$. For all $x\in \Omega$, all $t>0$ and all $r>0$,
\[
\left(\int_{r\leq \left\vert x-y\right\vert\leq 2r} \left\vert
\nabla_yK_t(y,x)\right\vert^2dy\right)^{1/2}\leq
cC_G t^{-\frac 12-\frac n4} \left(\frac r{\sqrt{t}}\right)^{\frac{n-2}2}
e^{-\beta\frac{r^2}t}.
\]
\item[$(b)$]
Assume that $\tau=1$. Then, for all $x\in \Omega$, all $0<t<1$ and all $r>0$,
\[
\left(\int_{r\leq \left\vert x-y\right\vert\leq 2r} \left\vert
\nabla_yK_t(y,x)\right\vert^2dy\right)^{1/2}\leq cC_G t^{-\frac 12-\frac
n4} \left(\frac r{\sqrt{t}}\right)^{\frac{n-2}2}
e^{-\beta\frac{r^2}t}.
\]
For all $x\in \Omega$, all $t>1$ and all $r>0$,
\[
\left(\int_{r\leq \left\vert x-y\right\vert\leq 2r} \left\vert
\nabla_yK_t(y,x)\right\vert^2dy\right)^{1/2}\leq cC_G t^{-\frac 12}
r^{\frac{n-2}2} e^{-\beta\frac{r^2}t}.
\]
\end{itemize}
\end{proposition}
As a consequence of Proposition \ref{gradient}, the following holds:
\begin{lemma} \label{heatint}
For all $x\in \Omega$, denote by $\delta(x)$ the distance from $x$
to $\partial \Omega$.
\begin{itemize} 
\item[$(a)$]
Under NBC, for all $x\in \Omega$,
\[
\int_{\Omega} \partial_tK_t(y,x)dy=0.
\]
\item[$(b)$]
Under DBC, for all $x\in \Omega$ for all $0<t<\tau$,
\[
\left\vert \int_{\Omega} \partial_tK_t(y,x)dy\right\vert \leq \frac Ct
e^{-\frac{\beta\delta(x)^2}{4t}}.
\]
Under DBC, if $\Omega$ is bounded and $\tau=1$, for all $x\in \Omega$
and all $t>1$,
\[
\left\vert \int_{\Omega} \partial_tK_t(y,x)dy\right\vert \leq \frac Ct.
\]
\end{itemize}
\end{lemma}
Under NBC, one has $e^{-tL}1=1$, whence assertion $(a)$ holds.

To prove the first part of assertion $(b)$, choose
$\psi_1\in C^{\infty}_0(\Omega)$ such that $\psi_1(z)=1$ if $d(z,y)\leq
\delta/4$,
$\psi_1(z)=0$ if
$d(z,y)\geq \delta/2$ and $\left\Vert \nabla \psi_1\right\Vert_{\infty}
\leq C/\delta$. Here $\delta=\delta(x)$. Define $\psi_2=1-\psi_1$.
Then,  one has
\[
\int_{\Omega} \partial_{t} K_t(y,x)dy=\int_{\Omega}
\partial_{t}K_t(y,x)\psi_1(y)dy+\int_{\Omega} \partial_{t} K_t(y,x)
\psi_2(y)dy.
\]
But Lemma \ref{gradient} shows that
\[
\begin{array}{lll}
\displaystyle \left\vert \int_{\Omega} L_yK_t(y,x)\psi_1(z)dz\right\vert
&=&\displaystyle \left\vert \int_{\Omega}
A\nabla_y K_t(y,x)\nabla_y\psi_1(y)dy\right\vert \\
\\
& \leq & \displaystyle C\int_{\frac{\delta}4 \leq d(z,y)\leq
\frac{\delta}2} \left\vert \nabla_y K_t(y,x)\right\vert
\left\vert \nabla_y\psi_1(y)\right\vert dy \\
\\
&\leq & \displaystyle Ct^{-\frac 12-\frac
n4}\left(\frac{\delta}{\sqrt{t}}\right)^{\frac{n-2}2}
e^{-\beta\frac{\delta^2}t} \delta^{\frac{n-2}2}\\
\\
&= &\displaystyle \frac {C}t
\left(\frac{\delta}{\sqrt{t}}\right)^{{n-2}}
e^{-\beta\frac{\delta^2}t}.
\end{array}
\]
Moreover,
\[
\begin{array}{lll}
\displaystyle \left\vert \int \partial_tK_t(x,y)\psi_2(y)dy\right\vert &
\leq & \displaystyle
\int_{\left\vert y-x\right\vert \geq \delta/2} \left\vert \partial_t
K_t(x,y)\right\vert dy \\
\\
& \leq & \displaystyle \frac Ct e^{-\beta\frac{\delta^2}t}.
\end{array}
\]
For the second part of assertion $(b)$, we have
$$
\int_\Omega |\partial_tK_t(x,y)|  \, dy \le \frac {C|\Omega|} t.
$$
\rule{2mm}{2mm} 

From these estimates and the subordination formula, we deduce the
following:
\begin{lemma} \label{poissint}
For all $x\in \Omega$, denote again by $\delta(x)$ the distance from $x$ to
$\partial \Omega$.
\begin{itemize}

\item[$(a)$]
Under NBC, for all $x\in \Omega$,
\[
\int_{\Omega} \partial_tp_t(x,y)dy=0.
\]
\item[$(b)$]
Under DBC, if $\tau=+\infty$, for all $x\in \Omega$,
\[
\left\vert \int_{\Omega} \partial_tp_t(x,y)dy\right\vert \leq \frac Ct
\left(1+\frac{ \delta(x)}{t}\right)^{-1}.
\]
Under DBC, if $\Omega$ is bounded and $\tau=1$, for all $x\in \Omega$
and $0<t<1$
\[
\left\vert \int_{\Omega} \partial_tp_t(x,y)dy\right\vert \leq \frac Ct
\left(1+\frac{ \delta(x)}{t}\right)^{-1}.
\]
\end{itemize}
\end{lemma}

Let us prove the second point of part $(b)$. By differentiating the
subordination formula, one has
$$
 \int_\Omega \partial_tp_t(x,y)\, dy=\frac
1{\sqrt{\pi}} \int_\Omega \int_0^{+\infty}
 \frac{2u}{t^2}\partial_sK_{s}(x,y){\mid_{s=\frac{t^2}{4u}}}
e^{-u}u^{-1/2}dudy.
$$
Break the integral at $u=t^2/4$. The part for $u\ge t^2/4$ is controlled
by $\frac Ct
\left(1+\frac{ \delta(x)}{t}\right)^{-1}$.  The part for $u\le t^2/4$
is bounded by $\frac c t \int_0^{t^2/4}e^{-u}u^{-1/2}du \le c $. Since $
t+\delta(x) \le 1 + {\rm diam}(\Omega)$, we obtain  
$c \le\frac Ct
\left(1+\frac{ \delta(x)}{t}\right)^{-1}$. This concludes the proof.
\rule{2mm}{2mm}  

We leave to the reader the care of studying what happens to 
regularity estimates for small time for $p_t(x,y)$ when $(G_1)$ holds.

 \section*{Appendix B: Elementary geometry of Lipschitz domains}
\label{geom}
\addcontentsline{toc}{section}{Appendix B:  Elementary geometry of
Lipschitz domains} 

A strongly Lipschitz domain is by definition a domain in $\reel^n$ whose
boundary is covered by a finite number of parts of Lipschitz graphs (up to
rotations) at most one them being infinite.  A special Lipschitz domain is
the domain above the  graph of a Lipschitz function defined on $\reel^{n-1}$.

Let $\Omega$ be a strongly Lipschitz domain.

\begin{enumerate}

\item There exists
a finite covering of $\reel^n$ by open sets
$U_1, U_2,
\ldots, U_s$ with at most one of them being infinite such that
for each $k$ either $U_k \cap \Omega=\emptyset$ or there is a special
Lipschitz domain
$\Omega_k$ and a rotation
$R_k$ in $\reel^n$ such that $U_k\cap \Omega= U_k\cap R_k(\Omega_k)$. 

\item  There exists a cube
$Q_0$ such that  either  $\Omega\subset Q_0$ or there is a rotation $R$ and a
special Lipschitz domain $\Omega_0$ such that $^cQ_0 \cap \Omega ={} ^cQ_0
\cap R(\Omega_0)$.

\item There are constants $\rho\in ]0,+\infty]$ and $C>0$ such that if $Q$ is
a cube centered in
$\Omega$ and $\ell(Q)\le\rho$ then $|Q\cap \Omega|\ge C|Q|$. When $\Omega$
is a unbounded, $\rho=\infty$.

\item There exists $\rho \in ]0,+\infty]$, such that  if
$Q$ is a type
$(b)$ cube and $\ell(Q)< \rho$,   there exists a cube
$\widetilde{Q}\subset
{}^c\Omega$ such that 
$\left\vert \widetilde{Q}\right\vert = \left\vert Q\right\vert$ and 
the distance from $\widetilde{Q}$ to $Q$ is comparable 
to the side length of $Q$. Furthermore, $\rho=\infty$ is ${}^c\Omega$ is
unbounded. 

\item Assume $\Omega$ is unbounded. Let
$Q$ be a cube  with $\ell(Q)\ge 1$, centered in $\Omega$ with $4Q\cap
\Omega\ne \emptyset$. There exists a  cube $Q'$ with $4Q'\subset
\Omega$,
$\ell(Q') =\ell(Q)$ and the distance between $Q$ and $Q'$ is comparable to
$\ell(Q)$.

\end{enumerate} 

The proof of 1. is classical and skipped. 

Point 2. follows easily: take
$Q_0$ as the smallest cube containing the bounded $U_k$'s in  point~1. 

To obtain $\rho=\infty$ in the proof of 3. when $\Omega$ is a special
Lipschitz domain is classical using ``vertical'' reflection. Localisation
gives us a finite $\rho$. To obtain $\rho=\infty$ when $\Omega$ is unbounded,
we argue as follows: let $Q_0$ be the cube of point 2. Let $Q$ be a cube
centered in $\Omega$ with $\ell(Q)>\rho$. 
If $\ell(Q) \le \lambda\ell(Q_0)$ for some $\lambda>1$ to be chosen, then for
$\tilde Q=
\frac
\rho{\ell(Q)} Q$  we have $ |Q\cap \Omega| \ge |\tilde Q\cap \Omega| \ge
C|\tilde Q| \ge C
\frac{\rho^n}{(\lambda\ell(Q_0))^n} |Q|$.
If $\ell(Q) \ge \lambda\ell(Q_0)$ and the center of $Q$ belongs to
$R(\Omega_0)$, then $|Q\cap \Omega| \ge |Q\cap \Omega\cap {}^cQ_0|
= |Q\cap R(\Omega_0)\cap {}^cQ_0| \ge |Q\cap R(\Omega_0)| - |Q_0|\ge
C|Q|- |Q_0|$ where $C$ is the constant obtained for the 
domain $R(\Omega_0)$. One chooses $\lambda$ so that $|Q_0| \le
C\rho^n/6^n$.  If $\ell(Q) \ge \lambda\ell(Q_0)$ and the center of $Q$
does not belong to
$R(\Omega_0)$, then this center belongs to $Q_0$ and one can find in $Q$ a
point in $R(\Omega_0)$ at distance less than $\ell(Q_0)$ from the center of
$Q$. It follows that $Q$ contains a cube of sidelength $\ell(Q)/3$ and
centered in $R(\Omega_0)$. We apply the above argument to that cube.
 
The proof of 4. is well-known if $\Omega$ is special Lipschitz or bounded. 
 See
\cite{stkrch}, p. 304.
By the same argument, one can see it holds for some $\rho$ finite for all
strongly Lipschitz domains. It remains to show that one can take drop the
finiteness of
$\rho$ if
$^c\Omega$ is unbounded.
In that case, let $Q$ be a type $(b)$ cube contained in $\Omega$ with
$\ell(Q)\ge
\rho$.  Pick $Q_0$, $R$ and $\Omega_0$ of point 2. In a basis $(e_1, \ldots ,
e_n)$,
$\Omega_0$ is $x_n\ge \varphi(x_1, \ldots, x_{n-1})$. We take
$\tilde Q =Q -c\ell(Q)R(e_n)$ for some appropriately chosen $c$ that depends
only on the domain $\Omega$. We leave to the reader the care of
verifying  that such a choice is possible. 

To see point 5. let
$Q$ is a cube of size greater than 1, centered in $\Omega$ with $4Q\cap
\Omega\ne \emptyset$. Arguing as above, we take $Q'=Q+c\ell(Q)R(e_n)$
where, since $\Omega$ is unbounded, one can pick $c$ large enough and
independent of
$Q$ such that
$Q'$ enjoys the desired properties.  Details are left to the reader.

\end{document}